\documentclass[12pt,a4]{article}

\usepackage{graphicx} 
\usepackage{amsmath, amsfonts, amsthm}
\usepackage[a4paper,left=2cm,right=2cm,top=2cm,bottom=2cm]{geometry}

\newtheorem{theorem}{Theorem}[section]
\newtheorem{prop}[theorem]{Proposition}
\newtheorem{lemma}[theorem]{Lemma}
\newtheorem{corollary}[theorem]{Corollary}

\theoremstyle{definition}
\newtheorem*{rem}{Remark}

\newcommand{\Ref}[1]{(\ref{#1})}

\newcommand{\X}{\mathcal{X}}
\newcommand{\Y}{\mathcal{Y}}
\newcommand{\Z}{\mathcal{Z}}
\newcommand{\F}{\mathcal{F}}

\newcommand{\nn}{\nonumber \\}

\usepackage{color}
\definecolor{IndianRed}{rgb}{0.8,0.36,0.36}
\definecolor{Green}{rgb}{0,0.8,0}
\definecolor{Blue}{rgb}{0,0,0.8}

\usepackage[
dvipdfm,%
colorlinks=true,%
urlcolor=IndianRed,%
linkcolor=Green,%
citecolor=IndianRed,%
pdfborder={0 0 1 [1]},
]{hyperref}

\numberwithin{equation}{section}

\begin{document}
\title{Partially directed paths in a wedge} 
\author{E J Janse van Rensburg$^1$\thanks{Please address
    correspondence to \texttt{rensburg@yorku.ca}}, T Prellberg$^2$ and
    A Rechnitzer$^{3,4}$.  \\[1ex]
  \footnotesize
  \begin{minipage}{13cm}
    $^1$ Department of Mathematics and Statistics,\\
    York University, 4700 Keele Street, Toronto, Ontario, M3J 1P3, Canada \\
    \texttt{rensburg@yorku.ca}\\[1ex]
    $^2$ School of Mathematical Sciences,\\
    Queen Mary College, Univ of London, Mile End Road, London E1 4NS,  UK\\
    \texttt{t.prellberg@qmul.ac.uk}\\[1ex]
    $^3$ Department of Mathematics and Statistics,\\
    The University of Melbourne, Victoria~3010, Australia.\\
    \texttt{a.rechnitzer@ms.unimelb.edu.au}\\  
    $^4$ Department of Mathematics\\
    University of British Columbia, Vancouver, BC~V6T~1Z2,
    Canada.\\
    \texttt{andrewr@math.ubc.ca}
  \end{minipage}
}

\maketitle  

\begin{center}
  \rule{10cm}{1pt}
\end{center}

\begin{abstract}
  The enumeration of lattice paths in wedges poses unique mathematical
  challenges.  These models are not translationally invariant, and the
  absence of this symmetry complicates both the derivation of a
  functional recurrence for the generating function, and solving for
  it.  In this paper we consider a model of partially directed walks
  from the origin in the square lattice confined to both a symmetric
  wedge defined by $Y = \pm pX$, and an asymmetric wedge defined by
  the lines $Y= pX$ and $Y=0$, where $p > 0$ is an integer.  We prove
  that the growth constant for all these models is equal to
  $1+\sqrt{2}$, independent of the angle of the wedge.  We derive
  functional recursions for both models, and obtain explicit
  expressions for the generating functions when $p=1$. From these we
  find asymptotic formulas for the number of partially directed
  paths of length $n$ in a wedge when $p=1$.

  The functional recurrences are solved by a variation of the kernel
  method, which we call the ``iterated kernel method''. This method
  appears to be similar to the obstinate kernel method used by
  Bousquet-M\'elou (see, for example, references \cite{BM02,BM03}).
  This method requires us to consider iterated compositions of the
  roots of the kernel.  These compositions turn out to be surprisingly
  tractable, and we are able to find simple explicit expressions for
  them. However, in spite of this, the generating functions turn out
  to be similar in form to Jacobi $\theta$-functions, and have natural
  boundaries on the unit circle.

  \vspace{2ex}
  
\end{abstract}

\begin{center}
  \rule{10cm}{1pt}
\end{center}


\section{Introduction}
\label{introduction}
The problem of counting random walks on lattices with various
restrictions is perhaps one of the oldest problems in enumerative
combinatorics, with a history that dates back at least 100 years
\cite{Andre1887}. It has also seen a great deal of recent activity,
particularly surrounding problems of counting random walks on the
slit-plane and quarter-plane \cite{Gessel1992, Niederhausen1998,
Fayolle1999, Barcucci2001, BM2001, BM2002a, BM02}. 

These models pose interesting mathematical problems, and powerful
methods have been developed in recent years to solve for the
generating functions of path problems.  Normally, these 
methods are a three step process:  First a recurrence is
determined, this is solved in the second step, and lastly,
the asymptotics of the number of paths are extracted. 

Perhaps the simplest and most studied model is Dyck paths.  While
there are numerous techniques for enumerating Dyck paths, the 
most powerful technique involves a recurrence for the generating
function which is solved and expanded to determine an explicit
expression.  If $d_n$ is the number of Dyck paths of half-length
$n$, then the generating function $g_t = \sum_{n\geq 0} d_n t^{n}$
satisfies the recurrence
\begin{equation}
g_t = 1 + tg_t^2
\end{equation}
with solution
\begin{equation}
g_t = \frac{2}{1+\sqrt{1-4t}} = \sum_{n=0}^\infty 
\left( {{2n}\atop{n}} \right) \frac{t^n}{n+1}
\end{equation}
so that $d_n$ is given by Catalan's number.

In this paper we follow a similar strategy to determine the 
generating function and asymptotic expressions for the number
of partially directed paths confined to a wedge.  The wedge
destroys translational invariance in the model, and both the
derivation of a recurrence for the generating function, and 
solving the generating function, poses difficult mathematical
problems.  Our strategy is in principle no different from the
above for Dyck paths - we shall derive functional recurrences
for the generating functions, solve those in special cases, and
then find the asymptotics for the number of paths.
Unfortunately, the problem for general wedges appears intractable,
and even in the cases that we do solve we encountered significant
difficulties.

Models of paths and walks frequently appear as simple models
of polymers in dilute solution in the physics literature
\cite{V98}.  The properties of polymers are in part determined
by their conformational entropy, and models of walks and paths
contributes to our understanding of the significance of 
the conformational entropy contributions in the free energy of
polymers.  These entropic contributions are important when
polymers are in confined geometries.  For example, the steric stabilisation of colloids by polymers results when polymers
are confined to the spaces between colloidal particles
\cite{N83}.  This situation have been modelled by studying
paths confined to the slab between two planes, see for example
\cite{DiMarzio1971, Brak2005}.

Lattice random walk models of polymers in confined geometries
are generally more tractable.  These models can generally be
solved, at least in principle, by a Bethe ansatz
or constant term formulation.  This technique have been used
to solve for random walks in a half-space and which interacts
with the boundary of the space \cite{MW40}.  Such random 
walk models, however, do not take into account the 
volume exclusion of monomers in a polymer.  A more 
realistic model is the self-avoiding walk \cite{MS93}.
This model is non-Markovian, and while much is known
about it from constructive \cite{K63,K64} and conformal 
invariance techniques (in two dimensions) \cite{Cardy1984},
solving it remains beyond the current techniques in combinatorics.

Self-avoiding walk models of polymers in confined spaces 
have not been solved (except in the most trivial of cases),
but there are some results in the literature.  For example, 
the exponential growth constant of self-avoiding walks  
in a wedge geometry is independent of the angle of the 
wedge \cite{HW85}. Additionally, for self-avoiding walks 
in wedges, conformal field theory have been used to
examine the dependence of scaling exponents on the wedge angle
\cite{Cardy1984, Duplantier1986}.

\begin{figure}[t]
\begin{center}
\includegraphics[width=\textwidth]{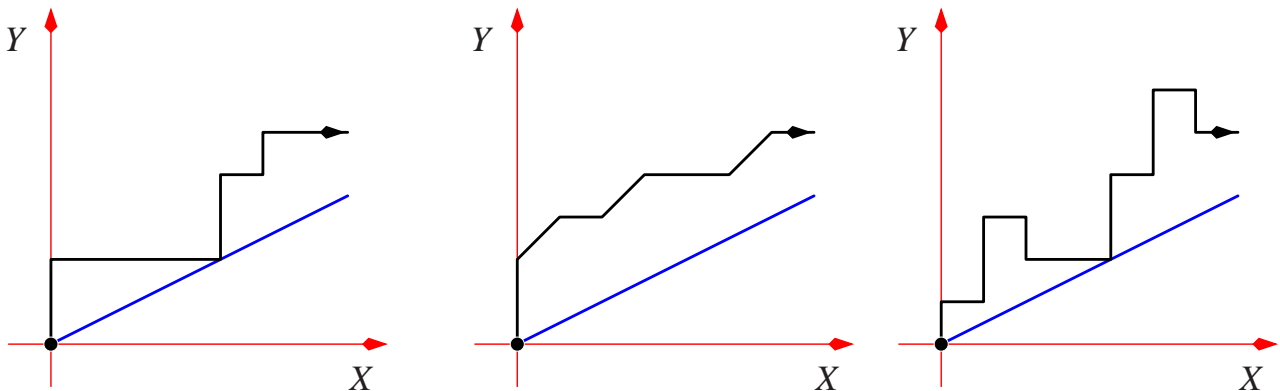}
\end{center}
\caption{Models of directed walks above a line $Y=+pX$. The walks are
  constrained to take only north and east steps (left); north, east
  and north-east steps (middle); and north, east and south steps (right).}
\label{fig1}
\end{figure}

The introduction of directedness in a self-avoiding walk in a wedge
may give models of directed or partially directed walks in a wedge
which are both self-avoiding and which may in some cases be solvable.
Models of directed and partially directed walks in wedge geometries
(see Figure~\ref{fig1}) have been studied previously \cite{D00,
  JvRL05, JvR05c}.  These models include directed paths in a wedge; a
model which is related to Dyck paths.  It is interesting that the
radius of convergence of the generating function is known in this model,
even for wedges with wedge-angles of irrational cotangent
\cite{JvRL05,JvR05c}.

In this paper we consider models of a partially directed path confined
in wedges (see Figure~\ref{fig2}).  These models are similar to the
directed path models in Figure~\ref{fig1}, however, they are also
substantially more challenging, since the path interacts with the
wedge on two sides, rather than on only one side.  As a result, it is
much harder to find their generating functions and analyse their
asymptotics.

\subsection{Directed and partially directed paths}

A directed walk on the square lattice is a path taking 
unit steps only in the north and east directions. Such objects 
are necessarily self-avoiding; they cannot revisit the same 
vertex. Partially directed paths may take unit steps only 
in the north, south and east directions with the further 
condition that no vertex is visited twice --- \emph{ie} they 
are self-avoiding.  Hence, north steps cannot be followed by 
south steps and \emph{vice-versa}. The generating function 
of such walks can be derived using standard techniques:
\begin{equation}
  W(t) = \sum_{n \geq 0} c_n t^n =  \frac{1+t}{1-2t-t^2},
\end{equation}
where $c_n$ is the number of walks of length $n$ and $t$ 
is the length generating variable.  An expansion of $W(t)$
of $W(t)$ in $t$ produces an explicit expression for $c_n$:
\begin{equation}
  \label{eqn num free walks}
  c_n  = \frac{1}{2}\left( (1+\sqrt{2})^{n+1} + (1-\sqrt{2})^{n+1}  \right).
\end{equation}
The exponential growth constant is the exponential rate at which
$c_n$ increases with $n$.  This is given by
\begin{equation}
  \mu = \lim_{n \to \infty} c_n^{1/n} = 1+\sqrt{2}.
\end{equation}
This is the most fundamental quantity in this model from a
statistical mechanics point of view.  The radius of convergence
of $W(t)$ is $\mu^{-1}$, and the limiting free energy
is the logarithm of the growth constant, $\kappa = \log \mu$,
this defines the explicit connection between the combinatorial
properties of the model and its thermodynamic properties.

We show some models of directed and partially directed paths in
wedge geometries in Figure~\ref{fig1}.  The model in
Figure~\ref{fig1}~(left) was considered in \cite{D00, JvR05a}. In
general the growth constant is a (non-trivial) function of the wedge
angle. The derivative of the free-energy with respect to the wedge
angle gives the moment of the entropic force exerted by the polymer on
the wedge and this was computed in \cite{JvR05a}. This model may be
also generalised by introducing an interaction between the line
$Y=+pX$ and the path, or by considering partially directed paths or
Motzkin paths instead \cite{JvR05b,JvR05c}.

\begin{figure}[t]
\begin{center}
  \includegraphics[width=\textwidth]{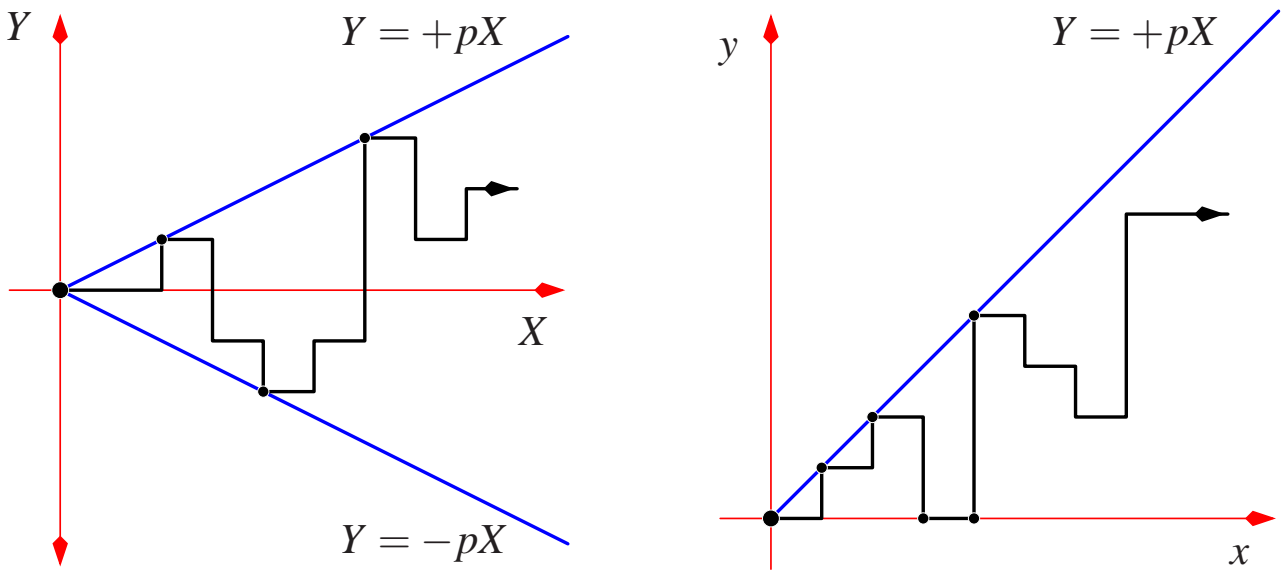}
\end{center}\caption{(left) The symmetric model of a partially directed 
path in a $p$-wedge formed by the lines $Y=\pm pX$.
(right) The asymmetric model of a partially directed path 
in a $p$-wedge formed by the line $Y=+pX$ and the $X$-axis.}
\label{fig2}
\end{figure}

In Figure~\ref{fig1}~(right) a partially directed path confined to the
wedge above the line $Y=pX$ and the $Y$-axis is proposed instead.
This model was considered in reference \cite{JvR05b}. If the partially
directed path is instead confined to the wedge between the $X$-axis
and the line $Y=pX$, then the model in Figure~\ref{fig2}~(right) is
obtained, which is the subject of this paper.

In particular, we consider the variants illustrated in
Figure~\ref{fig2} - firstly a model of a partially directed path in a
wedge formed by the lines $Y=\pm pX$ (we call this the {\it symmetric
  model} - see Figure~\ref{fig2}~(left)), and secondly a model of a
partially directed path in a wedge formed by the $X$-axis and the line
$Y=+pX$ (this is the {\it asymmetric model} - see
Figure~\ref{fig2}~(right)).

\begin{figure}[t]
\begin{center}
\includegraphics[width=6cm]{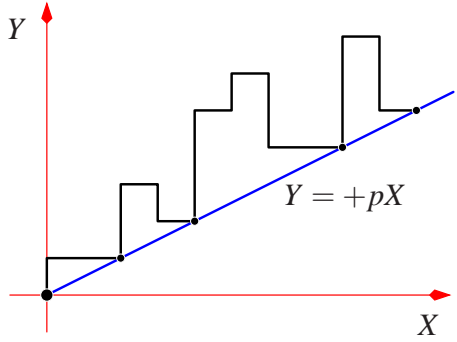}
\end{center}
\caption{A bargraph path above the line $Y=+pX$.  The generating
  function of this model is not known explicitly.  However, a set of
  equations derived in reference \cite{JvR05b} can be solved to
  determine the radius of convergence of the generating function for
  integer values of $p$.}
\label{fig3}
\end{figure}

The related model of a partially directed path in a wedge
with last vertex in the line $y=pX$ is illustrated
in Figure~\ref{fig3}. This is a {\it bargraph path} above 
the line $Y=+pX$.  This model was examined in reference 
\cite{JvR05b}, and while the generating function $g_p(t)$ 
is not known explicitly, it is given by
\begin{equation}
g_p(t) = \frac{h(t)}{1-t^2(1+h(t))}
\end{equation}
where $h(t)$ is an appropriate solution of the equation 
\begin{equation}
h(t) = t^{p+1}\left(1+h(t) \right)^p\left(1+\frac{h(t)}{1-t^2(1+h(t))}\right) ,
\end{equation}
where (as above), $t$ is conjugate to the number of edges 
in the path.  

The model in Figure~\ref{fig3} was also used as a model
of {\it adsorbing bargraphs} which interact with the 
line $Y=+pX$ \cite{D00}, and an asymptotic expression
for the adsorption critical point has been estimated 
in reference \cite{JvR05b} (the location of the singular 
point on the radius of convergence of the generating function).

\subsection{Partially Directed Paths in Wedges}

Consider the square lattice $\mathbb{Z}^2$ of points in 
the plane with integer coordinates.  Let $p>0$ be an integer.  
The symmetric $p$-wedge $\mathcal{V}_p$ is defined by
\begin{equation}
 \mathcal{V}_p = \left\{ (n,m)\in \mathbb{Z}^2 \, \big| \,
  \hbox{where $n\geq 0$ and $-pn \leq m \leq pn$} \right\}.
\end{equation}
The (asymmetric) {\it  $p$-wedge} is defined by 
\begin{equation}
  \mathcal{W}_p = 
  \left\{ (n,m)\in \mathbb{Z}^2 \,  \big| \, \hbox{where $n\geq 0$ and $0 \leq m \leq pn$} \right\}.
\end{equation}
Let $v_{p,n}$ (resp. $w_{p,n}$) be the number of partially 
directed walks in $\mathcal{V}_p$ (resp.  $\mathcal{W}_p$) 
of length $n$.

In the next section we establish some basic facts about 
the asymptotic growth of $v_{n,p}$ and $w_{n,p}$ as 
$n \to \infty$. In Section~\ref{section2} we find 
functional equations satisfied by the corresponding 
generating functions, which we solve in
Sections~\ref{sec sym} and~\ref{sec asym}.  We
then analyse the generating functions to determine
$v_{n,1}$ and $w_{n,1}$ to leading order.  We show 
in particular that
\begin{equation}
v_{n,1} = A_0 (1+\sqrt{2})^n + \frac{\sqrt{5}^n}{\sqrt{(n+1)^3}}
\left(A_1 + (-1)^n A_2 + O(1/n) \right)
\end{equation}
for the number of paths in a symmetric wedge when $p=1$ where
$A_0$, $A_1$ and $A_2$ are constants.
The asymmetric wedge poses more difficult mathematical 
problems, and we were only able to show that
\begin{equation}
w_{n,1} =  \frac{(1+\sqrt{2})^n}{\sqrt{n+1}} \left( B_0 + o(1)
\right),
\end{equation}
for some constant $B_0$.

\section{Growth constants in wedges} 

In this section we prove that the growth constant for 
partially directed walks is independent of the angle of 
the wedge and is equal to that of unrestricted partially
directed walks.  First, let $b_n$ be the number of 
partially directed walks in the wedge defined by the lines 
$X=0$ and $Y=0$, whose last vertex lies in the line $Y=0$. 
These paths are counted by the generating function 
$g_0(t)$ defined above, and singularity analysis gives 
the following lemma.

\begin{lemma}
  The growth constant of partially directed paths in the wedge defined
  by $X=0$ and $Y=0$ is
  \begin{equation}
    \lim_{n \to \infty} b_n^{1/n} = (1+\sqrt{2}) = \mu.
  \end{equation}
\end{lemma}

This result can be used to determine the growth constants of
partially directed walks in the wedges $\mathcal{V}_p$ and 
$\mathcal{W}_p$.   We first prove existence of the growth
constants.

\begin{lemma}
  For any given $p \in (0,\infty)$ the following limits exist:
  \begin{equation}
    \lim_{n \to \infty} v_{n,p}^{1/n} = \mu_p^v 
    \qquad \mbox{and} \qquad
    \lim_{n \to \infty} w_{n,p}^{1/n} = \mu_p^w.
  \end{equation}
  The limits satisfy
  \begin{equation}
    \mu_p^w \leq \mu_p^v \leq \mu =  (1+\sqrt(2))
  \end{equation}
\end{lemma}
\begin{proof}
  We have that $w_{n,p} \leq v_{n,p} \leq c_n$. Hence, if the 
above limits exist, we must have $\mu_p^w \leq \mu_p^v 
\leq \mu = (1+\sqrt{2})$.

  To show existence, we prove that the sequences are
super-multiplicative. Take any walk counted by $v_{n,p}$ 
and append a horizontal step, and any walk counted by 
$v_{m,p}$. This gives a walk of $n+m+1$ steps that lies within 
$\mathcal{V}_p$, and so is counted by $v_{n+m+1,p}$. Hence 
$v_{n,p} v_{m,p} \leq v_{n+m+1,p}$. A standard result 
(Fekete's lemma) on super-additive sequences (which we can apply 
by taking logarithms) then implies that $\mu_p^v$ exists. 
The proof for walks in $\mathcal{W}_p$ is identical.
\end{proof}

Next, we show that $\mu_p^w = \mu_p^v$, and we show that
they are equal to $1+\sqrt{2}$.

\begin{lemma}
  For any given $p \in (0,\infty)$ we have
  \begin{equation}
    b_{n}^N \leq w_{(\lceil np \rceil + nN + N), p}
  \end{equation}
  And hence $\lim_{n \to \infty} b_n^{1/n} \leq \mu_p^w$.
\end{lemma}
\begin{proof}
  Take any walk counted by $b_n$. By prepending $\lceil np \rceil
+1$ horizontal steps, this walk will fit inside the wedge, 
$\mathcal{W}_p$. Now append another horizontal step and a 
walk counted by $b_n$ --- repeat this until there are $N$ 
walks counted by $b_n$. This gives a walk counted by 
$w_{(\lceil np \rceil + nN + N), p}$. Thus we have the 
first inequality.  Taking logs and dividing by
$(\lceil np \rceil + nN + N)$ gives
  \begin{equation}
    \frac{N}{(\lceil np \rceil + nN + N)} \log b_{n} \leq
    \frac{1}{(\lceil np \rceil + nN + N)} \log w_{(\lceil np \rceil + nN + N), p} \\
  \end{equation}
Take the limit as $N \to \infty$ to obtain
  \begin{equation}
    \frac{1}{n} \log b_n \leq \log \mu_p^w.
  \end{equation}
Next, take the limit as $n \to \infty$ to complete the proof.
\end{proof}

By combining the above lemmas we can prove that the growth constant
for partially directed paths is independent of the wedge angle.

\begin{theorem}
  For any given $p \in (0,\infty)$ 
  \begin{equation}
    \mu_p^v = \mu_p^w = \mu = (1+\sqrt{2}).
  \end{equation}
\end{theorem}

This shows that the dominant asymptotic behaviour of the number of
walks is independent of the wedge angle. Below, we show that the
leading sub-dominant behaviour is also independent of the 
wedge angle (\emph{ie} for $p\geq1$).

\section{Functional equations for walks in wedges}
\label{section2}

\subsection{The symmetric wedge model}

Consider a model of partially directed paths in a symmetric wedge as
illustrated in Figure~\ref{fig2}~(right).  If $p$ is an integer or a
rational number, then the path may touch vertices in the lines $Y=\pm
pX$.  These vertices are {\it visits} in the lines $Y=\pm pX$. 
In the event that $p$ is an irrational number such visits 
cannot occur, however the path may approach arbitrarily close 
to the adsorbing lines (for large enough $X$-ordinate). 
In this paper we shall only consider the simplest version of 
this model, and we assume that $p$ is a positive integer.  
Even in this case the model is apparently intractable, and
we have only found the generating functions when $p=1$.

We will derive a functional equation satisfied by the generating
function of partially directed paths in $\mathcal{V}_p$ (those
illustrated in Figure~\ref{fig2}~(left)), by finding a recursive
construction, similar to those in \cite{BM96,BM02} (and elsewhere).

Let $x$ be the generating variable for horizontal edges in the 
path and let $y$ be the generating variable for vertical edges 
in the path. Introduce generating variables $a$ and $b$ to be
conjugate to the distances between the last vertex in the path and 
the line $Y=-pX$ and the line $Y=+pX$ respectively.
The generating function of the paths are now denoted by $g_p(a,b;x,y)
\equiv g_p(a,b)$ where the variables $x$ and $y$ are suppressed.

\begin{figure}[t!]
\begin{center}
\includegraphics[width=\textwidth]{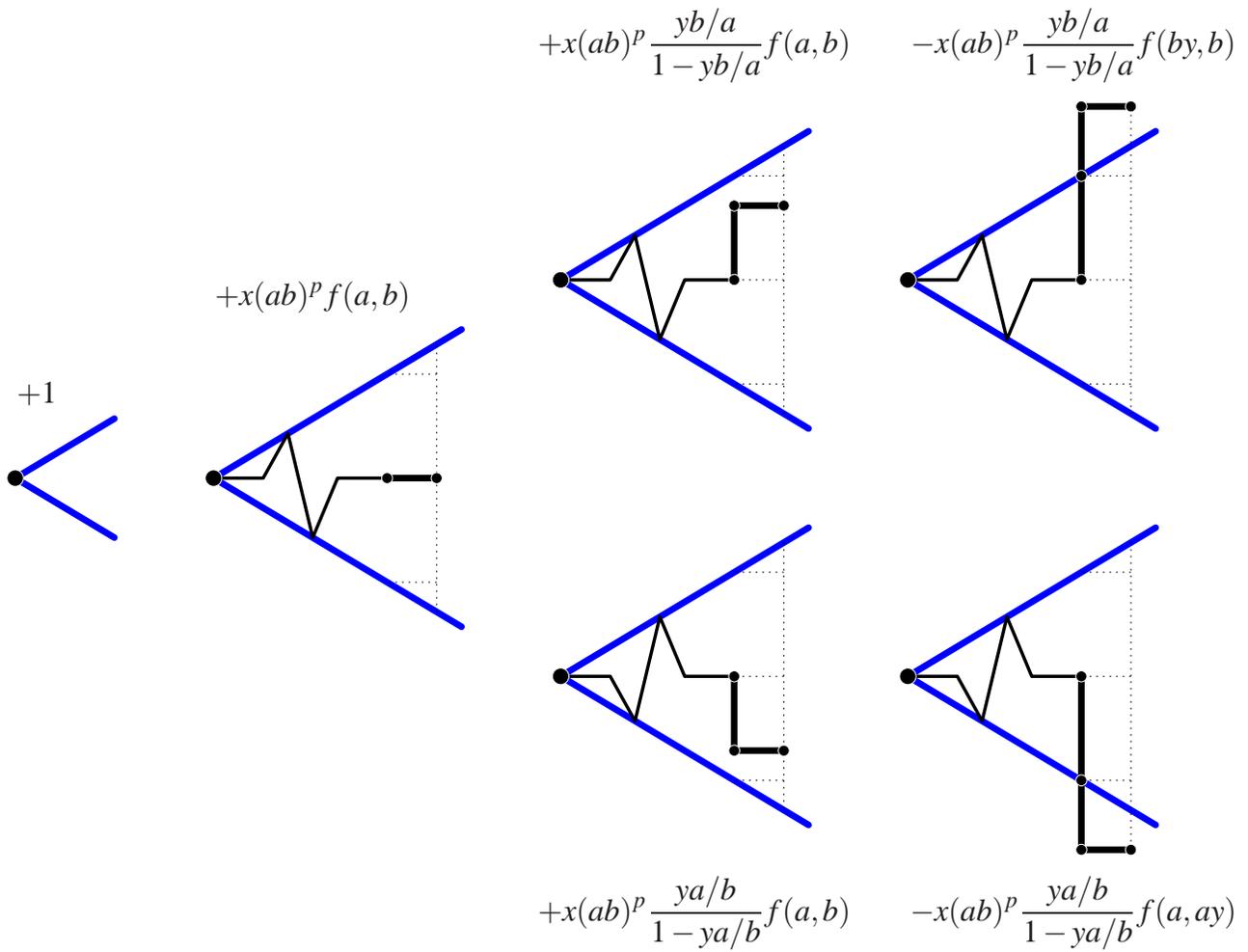}
\end{center}
\caption{Constructing partially directed walks in the wedge
  $\mathcal{V}_p$. Every walk is either a single vertex, or can be
  obtained from a shorter walk by appending a horizontal edge (left),
  or a run of north steps and a horizontal edge or a run of south
  steps and a horizontal edge (centre-top and -bottom). Care must be
  taken to not step outside the wedge when appending north or south
  steps (right-top and -bottom).}
\label{fig4}
\end{figure}

It turns out that the construction and resulting functional 
equation is simplified by considering only those partially 
directed walks that are either a single vertex (no edges) or 
end in a horizontal step. Let $f_p(a,b;x,y) \equiv f_p(a,b)$ 
be the generating function of such paths. It is simply related 
to $g_p(a,b)$ via
\begin{equation}
  \label{eqn gp}
  f_p(a,b) = 1 + x(ab)^p g_p(a,b).
\end{equation}
We now obtain a functional equation satisfied by $f_p$ by 
recursively constructing the paths column-by-column. Each path 
is either a single vertex, or can be constructed from a 
shorter path by appending either a horizontal step, or a 
sequence of up steps followed by a horizontal step, or a 
sequence of down steps followed by a horizontal step. See
Figure~\ref{fig4}.

Consider a path counted by $f_p(a,b)$, and see Figure~\ref{fig4}.
\begin{itemize}
\item Appending a single horizontal step to its end increases the
  distance of the end point from both wedge boundary lines by $p$.
  Hence the generating function of paths with a horizontal edge
  appended is $x (ab)^p f_p(a,b)$.
\item Appending an up step to the end of such a path increases
  the number of vertical steps by 1, increases the distance from the
  line $Y=-pX$ by 1 and decreases the distances from the line $Y=+pX$
  by 1. Hence such a path has generating function $y (b/a) f_p(a,b)$. 
  Hence appending some positive number of up steps gives 
  $\frac{yb/a}{1-yb/a} f_p(a,b)$. Appending a horizontal step 
  to the end of such a path gives (by the above reasoning) 
  $x(ab)^p\frac{yb/a}{1-yb/a}f_p(a,b)$.
\item Similarly appending some positive number of down steps followed
  by a horizontal step gives $x(ab)^p\frac{ya/b}{1-ya/b}f_p(a,b)$.
\end{itemize}

Unfortunately, when appending up or down steps it is possible 
that the resulting path will step outside of the wedge. Hence 
we must subtract off the contributions from such paths 
(Figure~\ref{fig4} right-top and -bottom).
\begin{itemize}
\item Consider a path that ends at a distance $h_+$ from the line
  $Y=+pX$. It we append more than $h_+$ up steps to the path then it 
  will leave the wedge. We can decompose the resulting path 
  into the original path with exactly $h_+$ up steps appended, 
  and an ``overhanging'' $\Gamma$ shaped path which is a 
  sequence of some positive number of up steps and a 
  horizontal step (see Figure~\ref{fig4} top-right).

  Appending exactly $h_+$ up steps to the path increases the distance
  from $Y=-pX$ by $h_+$, decreases the distance from $Y=-pX$ to zero.
  This gives the generating function $f_p(by,b)$.  The overhanging
  piece is (by the reasoning above) enumerated by $x(ab)^p
  \frac{yb/a}{1-yb/a}$.

  Hence the g.f. of walks that leave the wedge is given by $x(ab)^p
  \frac{yb/a}{1-yb/a} f(by,b)$.

\item Similarly when appending too many down steps we obtain
  configurations counted by $x(ab)^p \frac{ya/b}{1-ya/b} f_p(a,ay)$.
\end{itemize}

Using the above construction we arrive at the following theorem

\begin{prop}
  The generating function $f_p(a,b;x,y) \equiv f_p(a,b)$ of partially
  directed walks ending in a horizontal step in the wedge
  $\mathcal{V}_p$ satisfies the following functional equation:
  \begin{align}
    f_p (a,b) & =  1 + x(ab)^p f_p(a,b) \nn
    & + x(ab)^p \frac{yb/a}{1-yb/a} \left( f_p(a,b) - f_p(by,b) \right)\nn
    & + x(ab)^p \frac{ya/b}{1-ya/b} \left( f_p(a,b) - f_p(a,ay) \right)
    \label{eqn fp}
  \end{align}
  The generating function of all partially directed walks in
  $\mathcal{V}_p$ is given by
  \begin{equation}
    g_p(a,b) = x^{-1}(ab)^{-p} \left( f_p(a,b) - 1 \right)
  \end{equation}
\end{prop}

In the next section we turn to the problem of solving this 
functional equation.

\subsection{The asymmetric wedge model} 

Let us now turn to the construction of partially directed paths 
in the asymmetric wedge $\mathcal{W}_p$. Let the 
generating function of all partially directed walks in this 
wedge be denoted $k_p(a,b;x,y) \equiv k_p(a,b)$ where the 
variables $x$ and $y$ are suppressed. 

As above, the resulting functional equation satisfied by the
generating function is simpler if we consider only those walks 
that are either a single vertex or end in a horizontal step. 
Let this generating function be denoted $h_p(a,b;x,y)$. 
This is simply related back to $k_p$ by
\begin{equation}
  \label{eqn kp}
  h_p(a,b) = 1 + xa^p k_p(a,b).
\end{equation}

\begin{figure}[t!]
\begin{center}
\includegraphics[width=\textwidth]{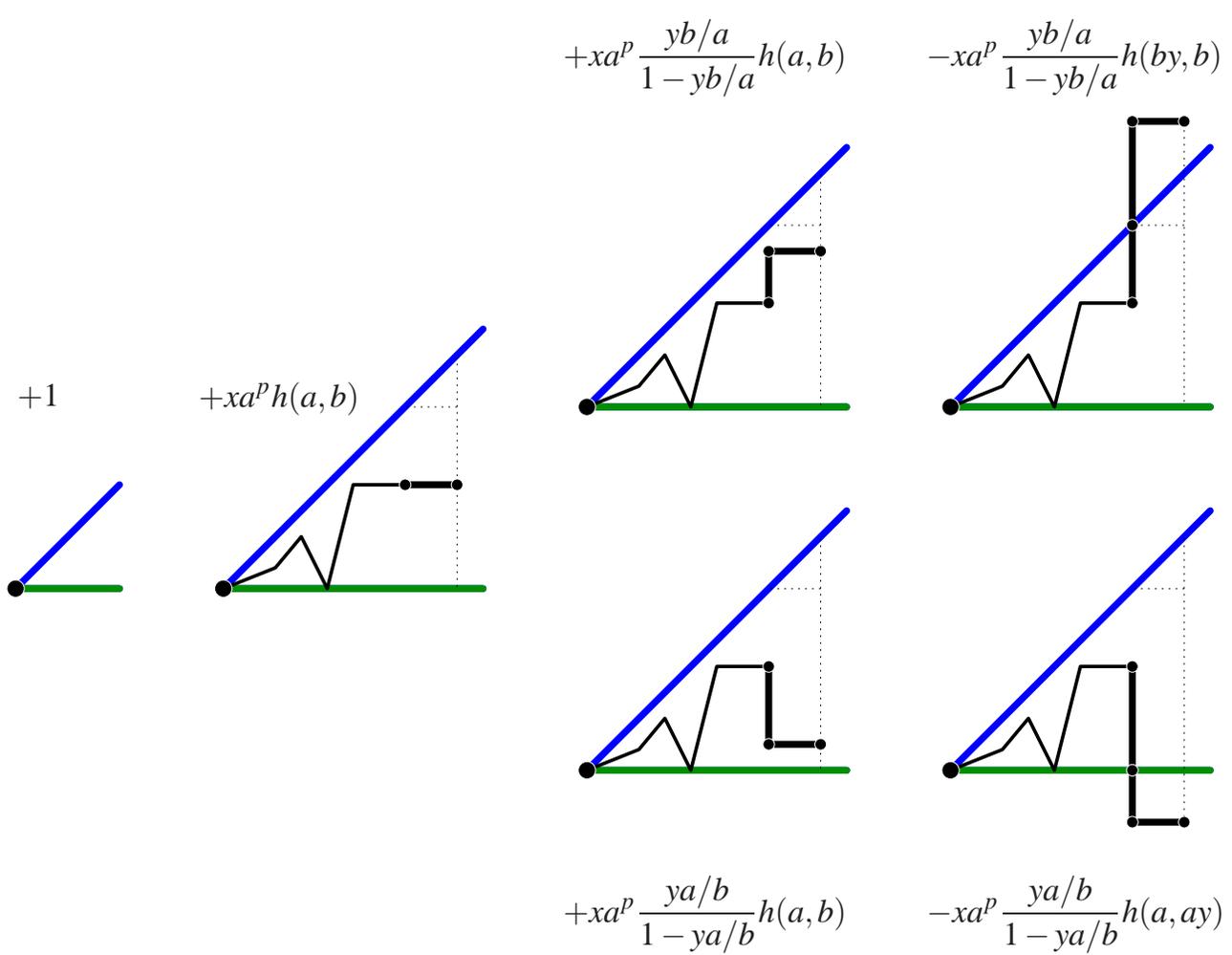}
\end{center}
\caption{Constructing partially directed walks in the asymmetric wedge
  $\mathcal{W}_p$. Every walk is either a single vertex, or can be
  obtained from a shorter walk by appending a horizontal edge (left),
  or a run of north steps and a horizontal edge or a run of south
  steps and a horizontal edge (centre-top and -bottom). Care must be
  taken to not step outside the wedge when appending north or south
  steps (right-top and -bottom).}
\label{fig5}
\end{figure}

We now use the same construction as was used above for the symmetric
case --- each walk is either a single vertex, or can be constructed
from a shorter walk by appending either a horizontal step, or a run of
up steps and a horizontal step, or a run of down steps and a
horizontal step --- see Figure~\ref{fig5}. Again care must be taken
not to step outside the wedge, and so those walks that do step outside
the wedge must be removed. Indeed the argument is \emph{almost}
identical to that used above, except that a horizontal step
contributes $x a^p$ instead of $x (ab)^p$, since a horizontal step
increases the distance from the line $Y=+pX$ by $p$, but does not
change the distance from the line $Y=0$.

The above construction gives the following theorem:

\begin{prop}
  \label{prop asym}
  The generating function $h_p(a,b;x,y) \equiv f_p(a,b)$ of partially
  directed walks ending in a horizontal step in the wedge
  $\mathcal{W}_p$ satisfies the following functional equation:
  \begin{align}
    h_p (a,b) & =  1 + xa^p h_p(a,b) \nn
    & + xa^p \frac{yb/a}{1-yb/a} \left( h_p(a,b) - h_p(by,b) \right)\nn
    & + xa^p \frac{ya/b}{1-ya/b} \left( h_p(a,b) - h_p(a,ay) \right)
    \label{eqn hp}
  \end{align}
  The generating function of all partially directed walks in
  $\mathcal{V}_p$ is given by
  \begin{equation}
    k_p(a,b) = x^{-1}a^{-p} \left( h_p(a,b) - 1 \right)
  \end{equation}
\end{prop}

We solve this equation in Section~\ref{sec asym}.

\section{Solving the symmetric case} 
\label{sec sym}
At first sight, one might try to solve equation~\Ref{eqn fp} by the
iteration method used in \cite{BM96}, however the coefficients of the
equation are singular when $a=by$ and $b=ay$. Multiplying 
both sides of the equation by $(a-by)(b-ay)$ gives a 
non-singular equation, however when we set $a=by$ or $b=ay$ 
the equation reduces to a tautology.

Instead we apply a variation of the kernel method, which we call 
{\it the iterated kernel method}. This appears to be similar 
in flavour to the ``obstinate kernel method'' used by 
Bousquet-M\'elou \cite{BM02,BM03}. We start by collecting all 
the $f(a,b)$ terms together on the left-hand side of the 
equation --- this gives the \emph{kernel form} of the equation:
\begin{equation}
K(a,b) \, f_p(a,b) = X(a,b) + Y (a,b) \, f_p(a,ya) + Z(a,b)\, f_p(yb,b),
\label{generic recursion}
\end{equation}
where the functions $K(a,b)$, $X(a,b)$, $Y(a,b)$ and $Z(a,b)$ are given by
\begin{subequations}
  \label{eqn coeff gen p}
  \begin{align}
    K(a,b) & = (b-ya)(a-yb)(1-x(ab)^p) - xy(ab)^p(a^2+b^2-2yab) , \\
    X(a,b) & =  (b-ya)(a-yb) \\
    Y(a,b) & = - xya^{p+1}b^p (a-yb) \\
    Z(a,b) & = - xya^pb^{p+1} (b-ya) 
  \end{align}
\end{subequations}
for each integer $p\geq 1$.  The function $K(a,b)$ is called the
\emph{kernel} of the equation. Note that the equation is symmetric
under interchange of $a$ and $b$:
\begin{equation}
  \label{eqn symms}
  f_p(a,b) = f_p(b,a) \qquad K(a,b) = K(b,a) \qquad X(a,b) = X(b,a)
  \qquad Y(a,b) = Z(b,a).
\end{equation}
We solve equation~\Ref{generic recursion} by substituting an infinite
number of pairs of $a$ and $b$ values that set the kernel $K(a,b)$ to
zero. While the method we describe below should work for general 
$p$, the resulting expressions are so complex that the 
process becomes intractable.  This will become clear even in the
case that $p=1$, which we simplified only after significant 
effort.

\subsection{Iterated kernel method for \texorpdfstring{${\cal V}_1$}{V1}}
\label{ssec it kern}

When $p=1$, the kernel becomes a quadratic function of 
$a$ and $b$ and we can explicitly write down (simple) 
expressions for its zeros. This is not generally true for 
larger values of $p$, and since simplifying our expressions
requires that compositions of the zeros of the kernel must 
simplify as well - the general $p$ case appears intractable
from a practical point of view.

Thus, we restrict ourselves to $p=1$. We write $f_1(a,b) \equiv
F(a,b)$ and the coefficients in equation~\Ref{generic recursion} become
\begin{subequations}
  \label{eqn coeff p1}
  \begin{align}
    \label{eqn kern p1}
    K(a,b) & = (xy^2a^2 - xa^2-y) b^2 + (1+y^2)ab - ya^2 \\
    X(a,b) & =  (b-ya)(a-yb) \\
    Y(a,b) & = - xya^2b (a-yb) \\
    Z(a,b) & = - xyab^2 (b-ya)
  \end{align}
\end{subequations}
Let $\beta_\pm(a;x,y) \equiv \beta_\pm(a)$ be the zeros of $K(a,b)$
with respect to $b$. Hence
\begin{equation}
  \label{eqn kern zero}
K(a,\beta_\pm (a)) = 0.
\end{equation}
Thus, setting $b=\beta_\pm (a)$ removes $F(a,b)$ from
equation~\Ref{generic recursion}.  This is the key idea behind the
``kernel method'' which has been used to solve equations of this type
(see \cite{B02} for example).

Unfortunately in this case, removing the kernel reduces 
the recurrence to an equation containing terms $F(a,ya)$ and 
$F(y\beta_\pm (a),\beta_\pm(a))$, which we cannot use immediately 
to solve for $F(a,b)$.  Similar situations have been studied 
before using the "obstinate kernel method" (\cite{BM02,BM03} 
for example).

The method we use appears to be similar to the obstinate kernel
method, except that instead of finding a finite number of pairs of
values of $a$ and $b$ to set the kernel to zero we must use an
infinite sequence of pairs. In this way, our ``iterated kernel
method'' is related both to the kernel method and perhaps also 
to the iterative scheme used in \cite{BM96}.

The roots $\beta_\pm (a)$ can be determined explicitly:
\begin{equation}
\beta_\pm (a) = \frac{a}{2} \left(
\frac{1+y^2\pm\sqrt{(1-y^2)(1-4xya^2-y^2)}}{y+xa^2-xy^2a^2}\right).
\label{eqn beta}
\end{equation}
Define the two roots:
\begin{align}
\beta_1 (a) \equiv \beta_-(a) & =  ya + O(xy^2a^3), \\
\beta_{-1} (a) \equiv \beta_+(a) & =  a/y + O(xy^{-2}a).
\end{align}
as power series in $a$. Later we require our solution to be a formal
power series in $t$ (after setting $x=y=t$) and one can confirm that
$\beta_1(a)$ defines a formal power series in $t$.

Since $a$ is a variable, we are able to substitute something else for
it; substituting $a \mapsto \beta_1(a)$ into 
equation~\Ref{eqn kern zero} gives
\begin{equation}
  K(\beta_1(a),\beta_1(\beta_1(a)) = 0.
\end{equation}
Hence the pair $(a,b) 
= \big(\beta_1(a), \beta_1(\beta_1(a))\big)$ also sets 
the kernel to zero. We can continue in this way. Hence we 
need to define the repeated composition of $\beta_1(a)$ with itself:
\begin{equation}
  \beta_n (a) = \beta_1^{(n)} (a) = 
  \underbrace{\beta_1 \circ \beta_1 \circ \ldots \circ\beta_1}_{n} (a).
\end{equation}
Note that 
\begin{align}
  \beta_{-1}\circ\beta_1 (a) & = \beta_1 \circ \beta_{-1} (a) = a,
  \qquad \mbox{and}\\
  \beta_n(a) & =  ay^n + O(xy^{n+1} a^3).
\end{align}
There is no finite value of $n$ such that $\beta_n = \beta_0$. If we
further define $\beta_0 (a) = a$ and $\beta_{-n} (a)$ by composition
of $\beta_{-1} (a)$, then the functions $\{\beta_n \, | \, n\in
\mathbb{Z} \}$ form an infinite group with identity $\beta_0 $
and inverses $\beta_n \circ \beta_{-n} = \beta_0$.

These observations are enough to iterate the functional equation to
find a solution.  Set $b = \beta_1 (a)$ in~equation \Ref{generic
  recursion}, and set $a = \beta_n (a)$ for any finite $n\geq 0$.
Then since $K(\beta_n(a), \beta_{n+1}(a)) = 0$, we have
\begin{align}
  X(\beta_n (a),\beta_{n+1} (a)) 
  & +  Y (\beta_n (a),\beta_{n+1} (a))  
  \, F(\beta_n (a),y\beta_n (a)) \nn
  & +  Z (\beta_n (a),\beta_{n+1} (a)) 
  \, F(y\beta_{n+1} (a),\beta_{n+1} (a)) = 0 .
\end{align}
We can then solve this equation for $F(\beta_n (a),y\beta_n (a))$:
\begin{align}
  F(\beta_n (a),y\beta_n (a)) =
  & - \left[ \frac{X(\beta_n (a),\beta_{n+1} (a))}{
      Y (\beta_n (a),\beta_{n+1} (a))} \right] \nonumber \\
  & - \left[\frac{Z (\beta_n (a),\beta_{n+1} (a))}{
      Y (\beta_n (a),\beta_{n+1} (a))} \right]
  F(y\beta_{n+1} (a),\beta_{n+1} (a))
  \label{eqn 21} 
\end{align}
We can simplify the above by defining
\begin{eqnarray}
\F_n (a) & = & F(\beta_n (a),y\beta_n (a)) 
= F(y\beta_n (a),\beta_n (a)), \nonumber \\
\X_n (a) & = & - \left[ \frac{X(\beta_n (a),\beta_{n+1} (a))}{
   Y (\beta_n (a),\beta_{n+1} (a))} \right],  \\
\Z_n (a) & = & - \left[\frac{Z (\beta_n (a),\beta_{n+1} (a))}{
   Y (\beta_n (a),\beta_{n+1} (a))} \right],
\end{eqnarray}
where we have made use of the symmetry $F(a,b) = F(b,a)$. While this
symmetry is not essential, it does make the solution substantially
simpler. Instead of exploiting this symmetry we could iterate again to
find $F(\beta_n (a),y\beta_n (a))$ in terms of $F(\beta_{n+2}
(a),y\beta_{n+2} (a))$. Indeed this is what is required to solve walks
in the asymmetric wedge $\mathcal{W}_1$ (see Section~\ref{sec asym} below).

Equation \Ref{eqn 21} may be written as
\begin{equation}
\F_n (a) = \X_n (a) + \Z_n (a) \F_{n+1} (a) .
\end{equation}
Starting at $n=0$, this can be iterated to get a series
solution for $\F_0 (a)$:
\begin{equation}
F(a,ya) = \F_0 (a) = 
\sum_{n=0}^\infty \X_n (a) \prod_{k=0}^{n-1} \Z_k (a),
\end{equation}
where we have \emph{assumed} that the above sum converges (we will
show that this is the case). This also gives $F(yb,b)$:
\begin{equation}
F(yb,b) = F(b,yb) = \F_0 (b) =
\sum_{n=0}^\infty \X_n (b) \prod_{k=0}^{n-1} \Z_k (b) .
\end{equation}
This allows us to write down the solution for $F(a,b)$:
\begin{equation}
  f_p (a,b) = 
  \frac{X(a,b)}{K(a,b)} +
  \frac{Y(a,b)}{K(a,b)}\sum_{n=0}^\infty \X_n (a) 
  \prod_{k=0}^{n-1} \Z_k (a) +
  \frac{Z(a,b)}{K(a,b)}\sum_{n=0}^\infty \X_n (b) 
  \prod_{k=0}^{n-1} \Z_k (b) .
\end{equation}

Of course, the above ``solution'' still contains many complicated
algebraic functions in the form of the $\beta_n(a)$. It is quite
surprising (at least to the authors!) is that these functions can be
drastically simplified.

\subsection{An explicit expression for \texorpdfstring{$f_1(1,1)$}{f1(1,1)}}
\label{ssec explicit sym}

We have outlined above the iterated kernel method that we shall 
use to write down the generating function $f_1(a,b) = F(a,b)$. 
We are primarily interested in the number of paths (and not 
the location of their endpoints), so we will actually focus 
on the function $F(1,1)$.

We start by considering the $\beta_n (a)$ functions. It is quite
surprising that while $\beta_n(a)$ is (upon superficial inspection for
small $n$) very complicated, its reciprocal appears relatively simple.
Examining equations~\Ref{eqn kern p1} and~\Ref{eqn beta} one obtains
\begin{equation}
\frac{1}{\beta_1(a)} + \frac{1}{\beta_{-1}(a)}
= \frac{1+y^2}{y} \frac{1}{a} . 
\end{equation}
Substituting $a = \beta_{n-1} (a)$ in the above, and using
the group properties of $\beta_n$ leads to the following
three term recurrence for $\beta_n$:
\begin{equation}
\frac{1}{\beta_n} = \frac{1+y^2}{y} \frac{1}{\beta_{n-1}} 
- \frac{1}{\beta_{n-2}}.
\end{equation}
Since $\beta_0$ is the identity, and $\beta_1$ is given
explicitly by $\beta_-(a)$ in equation \Ref{eqn beta}, the
recurrence above can be iterated to get a solution for
$\beta_n (a)$:
\begin{equation}
  \label{eqn betan}
  \frac{1}{\beta_n (a)} = 
  \frac{y(1-y^{2n})}{y^n(1-y^2)}\frac{1}{\beta_1 (a)} -
  \frac{y^2(1-y^{2n-2})}{y^n(1-y^2)}\frac{1}{a} .
\end{equation}
By using the expressions for $X(a,b)$, $Y(a,b)$ and $Z(a,b)$
in equation~\Ref{eqn coeff p1} to determine $\X(a,b)$ and $\Z(a,b)$,
one obtains
\begin{equation}
F (a,ya) = \sum_{n=0}^\infty (-1)^n
\left[ \frac{\beta_{n+1} - y \beta_n}{xya\beta_n \beta_{n+1}} \right]
\prod_{k=0}^{n-1} \left( 
\frac{\beta_{k+1} - y \beta_k}{\beta_k - y \beta_{k+1}}
\right) .
\end{equation}
Substituting the expression for $\beta_n(a)$ given in
equation~\Ref{eqn betan} and simplifying gives:
\begin{align}
  \frac{\beta_{n+1} - y \beta_n}{xya\beta_n \beta_{n+1}}
  & = y^n \left[ \frac{1}{a} - \frac{y}{\beta_1} \right] , \qquad
  \mbox{and} \\
  \frac{\beta_{k+1} - y \beta_k}{\beta_k - y \beta_{k+1}}
  &= y^{2k+1}
  \left[ \frac{1}{xya^2} - \frac{1}{xa\beta_1} - 1 \right] .
\end{align}
Using these, one can get an explicit expression for the generating
function $F(a,ya)$:
\begin{equation}
F (a,ya) = \left[\frac{1}{xya^2} - \frac{1}{xa\beta_1} \right]
\sum_{n=0}^\infty (-1)^n y^{n(n+1)} 
\left( \frac{1}{xya^2} - \frac{1}{xa\beta_1} - 1 \right)^n.
\end{equation}
By defining 
\begin{equation}
  \label{eqn qdef}
  Q(a;x,y) = \left( \frac{1}{xa^2} - \frac{y}{xa\beta_1} - y \right).
\end{equation}
the above expression for $F(a,ya)$ can be further simplified to 
\begin{equation}
  \label{eqn faya}
  F (a,ya) = \left[1 + \frac{Q(a;x,y)}{y} \right]
  \sum_{n=0}^\infty (-1)^n y^{n^2} Q(a;x,y)^n.
\end{equation}
Using the $a \leftrightarrow b$ symmetry of $F(a,b)$, we can get a similar
expression for $F(yb,b)$, and so finally $F(a,b)$.
\begin{multline}
  f(a,b)  = \frac{X(a,b)}{K(a,b)}
  + \frac{Y(a,b)}{K(a,b)} \left(1 + \frac{Q(a)}{y} \right) \sum_{n \geq 0} (-1)^n Q(a)^n y^{n^2} \\
   + \frac{Z(a,b)}{K(a,b)} \left(1 + \frac{Q(b)}{y} \right) \sum_{n \geq 0} (-1)^n Q(b)^n y^{n^2}
\end{multline}

We can reduce the above equation by considering only the number of
walks of length $n$ (by setting $a=b=1, x=y=t$):
\begin{prop}
  The generating function of partially directed walks ending in a
  horizontal step in the wedge $\mathcal{V}_1$ is
  \begin{equation}
    f_1 (1,1) = \frac{1-t}{1-2t-t^2} 
    -  \frac{1-t^2-\sqrt{(1-t^2)(1-5t^2)}}{1-2t-t^2}
    \sum_{n=0}^\infty (-1)^n t^{n^2} Q(1;t,t)^n ,
  \end{equation}
  where $t$ counts the number of edges and 
  \begin{equation}
    Q(1;t,t) = (1-3t^2-\sqrt{(1-t^2)(1-5t^2)})/2t .
  \end{equation}
  The generating function of all paths in $\mathcal{V}_1$ is then
  found using equation~\Ref{eqn gp}:
  \begin{equation}
    g_1 (1,1) = \frac{1+t}{1-2t-t^2} 
    -  \frac{1-t^2-\sqrt{(1-t^2)(1-5t^2)}}{t(1-2t-t^2)}
    \sum_{n=0}^\infty (-1)^n t^{n^2} Q(1;t,t)^n.
    \label{eqn 36}
  \end{equation}
\end{prop}

Firstly we note that $F(a,ya)$ counts all partially directed paths in
the wedge $\mathcal{V}_1$ whose last vertex ends in the line $Y=-pX$.
Additionally we note that the generating function $Q(a;x,y)/y$ counts
the number of partially directed paths starting at the origin, lying
on or above the line $Y=-pX$ and whose last vertex lies in the line
$Y=-pX$. Hence $Q(a)/y$ counts a very similar set of paths to 
$F(a,ya)$, except that the paths counted by $Q$ are not confined 
by the line $Y=pX$.

In light of the above interpretation of the function $Q$, we 
expended considerable effort to uncover a more direct 
combinatorial derivation of the alternating sum 
in equation~\Ref{eqn faya}. There appears to be some 
inclusion-exclusion process underlying this, but unfortunately
we have not made progress in this respect.

\subsection{Asymptotics for \texorpdfstring{$p=1$}{p=1}}

The asymptotics of the number of partially directed paths in the
symmetric wedge with $p=1$ can be analysed by examining the
singularities of the generating function $g_1(1,1)$ in equation
\Ref{eqn 36}.  Singularities arise either as zeros of the factor
$(1-2t-t^2)$ in equation~\Ref{eqn 36}, or as singularities in
$\sqrt{(1-t^2)(1-5t^2)}$, or as singularities in the series
$\sum_{n=0}^\infty (-1)^n t^{n^2} Q(1;t,t)^n$.

An examination of $g_1(1,1)$ shows that it has simple poles at 
the solution of $(1-2t-t^2)=0$, or when $t = -1 \pm\sqrt{2}$. 
We note that $\sqrt{(1-t^2)(1-5t^2)}$ has branch-points 
(square root singularities) at $t = \pm 1$ and again at 
$t=\pm 1/\sqrt{5}$.  The series $\sum_{n=0}^\infty (-1)^n 
t^{n^2} Q(1;t,t)^n$ is a Jacobi $\theta$-function and it 
is convergent inside the unit circle except at singularities of
$Q(1;t,t)$; that is, when $t=\pm 1/\sqrt{5}$.

The dominant singularity is the simple pole at $\sqrt{2}-1$, 
while the next sub-dominant contributions to the asymptotics 
will be given by the singularities at $t=\pm 1/\sqrt{5}$. 
These two sub-dominant singularities will give a parity effect. 
The contributions from these singularities allow us to write 
down the asymptotic form of $v_{n,1}$.

\begin{prop}
  The number of paths in the wedge $\mathcal{V}_1$ is asymptotic to
  \begin{equation}
    v_{n,1} = A_0 \left(1+\sqrt{2}\right)^n
    + \frac{5^{n/2}}{(n+1)^{3/2}} 
    \Big(A_1 + (-1)^n A_2 + O(1/n)\Big) .
    \label{eqn 37}
  \end{equation}
  Where the constants are
  \begin{subequations}
    \begin{align}
      A_0 & = 0.27730985348603118827\ldots , \\
      A_1 & =  3.71410486533662324953\ldots, \\
      A_2 & =  0.20697997020804157910\ldots.
    \end{align}
  \end{subequations}
\end{prop}

We note that the constants were derived by expanding the expression
for $g_1(1,1)$ about $t=\sqrt{2}-1$ and $t=\pm 1/\sqrt{5}$ (or rather
the first 40 or so terms of the sum). These were then checked using
both Bruno Salvy's \emph{gdev} package for Maple \cite{Salvy91b} and
by direct examination of $v_{n,1}$ for $n \leq 1000$. The above
formula is quite precise and it correctly estimates $v_{10,1},
v_{20,1}, v_{30,1}$ and $v_{40,1}$ to within $7\%, 1\%, 0.2\%$ and
$0.06\%$ respectively.

Note that the above result implies that walks in the wedge
$\mathcal{V}_p$ have the same dominant asymptotic behaviour as walks
with no bounding wedge (see equation~\Ref{eqn num free walks}).
Since the number of walks in any wedge $\mathcal{V}_p$ for
$1 \leq p < \infty$ is bounded between the number of walks in
$\mathcal{V}_1$ and partially directed walks with no bounding wedge,
we have the following result:
\begin{corollary}
  \label{cor symm bound}
  The number of partially directed walks in the wedge $\mathcal{V}_p$,
  $c_n^{(p)}$ obeys the following inequality
  \begin{equation}
    0.2773\ldots \leq \lim_{n \to \infty} \frac{c_n^{(p)}}{(1+\sqrt{2})^n} \leq
    (1+\sqrt{2})/2 = 1.2071\ldots
  \end{equation}
  for any $1 \leq p < \infty$.
\end{corollary}


\section{Partially Directed Paths in the Asymmetric Wedge}
\label{sec asym}

In this section we turn our attention to the model in
Figure~\ref{fig2}~(right).  The partially directed path is 
confined to an asymmetric wedge, $\mathcal{W}_p$, and 
its generating function does not have the 
$a \leftrightarrow b$ symmetry we have exploited in solving
for $f_1(a,b)$ in the previous section.

We proceed by examining the generating function of walks that end
in a horizontal step. The functional equation for these walks 
is given in Proposition~\ref{prop asym} and we can arrange 
equation \Ref{eqn hp} in kernel form:
\begin{equation}
K(a,b) \, h_p(a,b) = X(a,b) + Y (a,b) \, h_p(a,ya) + Z(a,b)\, h_p(yb,b)
\label{eqn Hkern}
\end{equation}
where
\begin{subequations}
  \label{eqn XXh} 
  \begin{align}
    K(a,b) & = (b-ya)(a-yb)(1-xa^p)-xya^p(a^2+b^2-2yab) , \\
    X(a,b) & = (b-ya)(a-yb), \\
    Y(a,b) & =  - xya^{p+1}(a-yb), \\
    Z(a,b) & =  - xya^pb(b-ya) .
  \end{align}
\end{subequations}
This functional equation is very similar to that of the symmetric
wedge given in equation~\Ref{generic recursion}. However we no longer
have $a \leftrightarrow b$ symmetry and this means that we have to
work quite a bit harder and we concentrate only on the case $p=1$. We
will write $h_1(a,b) \equiv H(a,b)$ for the remainder of this section.

\subsection{Solving for \texorpdfstring{$H(a,b)$ when $p=1$}{H(a,b)
    when p=1}}

For the remainder of this section we concentrate on the case $p=1$ and
walks in the $45^\circ$ wedge $\mathcal{W}_1$. The kernel $K(a,b)$
(given in equation~\Ref{eqn XXh}) is no longer symmetric in $a$ and
$b$, nor is the desired generating function $H(a,b)$. In order to
repeat the iterated kernel method as described in 
Section~\ref{ssec it kern} we must now consider the solutions 
of the kernel as functions of $a$ and $b$. These solutions 
are defined by $K(a, \beta(a)) = 0$ and
$K(\alpha(b),b) = 0$:
\begin{equation}
\hspace{-2cm}
\beta_\pm (a) = \frac{a}{2y} \left[
1+y^2-x(1-y^2)a \pm \sqrt{(1-y^2)((1-xa)^2-y^2(1+xa)^2)}\right]
\end{equation}
and
\begin{equation}
\alpha_\pm (b) = \frac{b}{2} \left[
\frac{1+y^2\pm\sqrt{(1-y^2)(1-y^2-4xyb)}}{y+x(1-y^2)b}
\right].
\end{equation}
One can confirm that $\alpha_-(b)$ and $\beta_-(a)$ both define 
formal power series in $b$ and $a$ (respectively). Additionally 
these same choices (when $x=y=t$) also define formal power series 
in $t$ --- which we will require for our solution.  Write these 
as $\alpha_1(b)$ and $\beta_1(b)$, and the other roots as 
$\alpha_{-1}(b)$ and $\beta_{-1}(a)$.

In Section~\ref{ssec it kern} we considered composing $\beta(a)$ 
with itself, however due to the asymmetry of the kernel we now 
need to consider mixed compositions $\beta(\alpha(b))$ and 
$\alpha(\beta(a))$.  Indeed, we find that
\begin{subequations}
  \label{eqn alpbet comp}
  \begin{align}
    \alpha_{\pm1} ( \beta_{\mp1} (a) ) & =  a , \\
    \beta_{\pm1} ( \alpha_{\mp1} (b)) & = b.
  \end{align}
\end{subequations}
We will need the function $\gamma(a) = \alpha_1( \beta_1 (a))$, 
and define its nested composition by $\gamma_n(a) =
\gamma(\gamma_{n-1}(a))$ with $\gamma_0(a) = a$. Note that
\begin{equation}
  \gamma_n(a) = y^{2n} a + O( x y^{2n} a^2).
\end{equation}

We can now repeat the iterated kernel method in the new asymmetric
setting. Setting $b=\beta_1(a)$ in equation~\Ref{eqn Hkern} gives
\begin{equation}
  \label{eqn habet}
  0 = X(a,\beta_1(a)) + Y(a,\beta_1(a)) H(a,ya) + Z(a,\beta_1(a))
  H(y\beta_1(a), \beta_1(a)).
\end{equation}
Since there is apparently not a simple relation between
$H(by,b)$ and $H(b,by)$, this equation cannot be iterated
to find a solution.  Instead, it turns out that the other
roots of the kernel must be considered as well.

Setting $a = \alpha_1(b)$ gives:
\begin{equation}
  \label{eqn halpb}
  0 = X(\alpha_1(b),b) + Y(\alpha_1(b),b) H(\alpha_1(b),y\alpha_1(b)) 
  + Z(\alpha_1(b),b) H(yb,b).
\end{equation}
Now set $b = \beta_1(a)$ in the above equation
\begin{multline}
  \label{eqn hgambet}
  0 = X(\gamma(a),\beta_1(a)) + Y(\gamma_1(a),\beta_1(a)) H(\gamma_1(a),y\gamma_1(a))
  \\
+ Z(\gamma_1(a),\beta_1(a) ) H(y\beta_1(a),\beta_1(a)).
\end{multline}
We can now eliminate $H(y\beta_1(a),\beta_1(a))$ between
equations~\Ref{eqn habet} and~\Ref{eqn hgambet} and solve for
$H(a,ya)$. This gives
\begin{multline}
  \label{eqn hayarec}
  H(a,ya) =
  - \left[\frac{X(a,\beta_1(a))}{Y(a,\beta_1(a))}\right]
  + \left[ \frac{Z(a,\beta_1(a))}{Y(a,\beta_1(a))}\right]
  \left[ \frac{X(\gamma_1 (a),\beta_1 (a))}
    {Z(\gamma_1(a),\beta_1 (a))}\right]\\
  + \left[ \frac{Z(a,\beta_1(a))}{Y(a,\beta_1(a))}\right]
  \left[ \frac{Y(\gamma_1(a),\beta_1 (a))}{
      Z(\gamma_1(a),\beta_1 (a))}\right] 
  H(\gamma_1 (a),y\gamma_1 (a)) .
\end{multline}
We can now iterate the above equation by substituting $a =
\gamma_{n-1}(a)$. Define
\begin{subequations}
  \begin{alignat}{2} 
    {\X}_n(a) & = 
    - \frac{X(\gamma_n(a),\beta_1(\gamma_n(a)))}
      {Y(\gamma_n(a),\beta_1(\gamma_n(a)))}
      && = \frac{\beta(\gamma_n) - y \gamma_n}{xy \gamma_n^2},
      \\
    {\Y}_n(a) & = 
    \frac{Z(\gamma_n(a),\beta_1(\gamma_n(a)))}
      {Y(\gamma_n(a),\beta_1(\gamma_n(a)))}
      && = \frac{\beta(\gamma_n)}{\gamma_n} \left(
        \frac{\beta(\gamma_n)-y\gamma_n}{\gamma_n - y \beta(\gamma_n)}
      \right) ,\\
    {\Z}_n(a) & = 
    \frac{X(\gamma_{n+1}(a),\beta_1(\gamma_n(a)))}{
        Z(\gamma_{n+1} (a),\beta_1 (\gamma_n(a)))}
      && = - \frac{\gamma_{n+1} - y \beta(\gamma_n)}{x y \gamma_{n+1} \beta(\gamma_n)}
 ,\\    
    \mathcal{A}_n(a) & =  
    \frac{Y(\gamma_{n+1} (a),\beta_1
      (\gamma_n(a)))}{ Z(\gamma_{n+1} (a)),\beta_1
      (\gamma_n(a)))}
    && = \frac{\gamma_{n+1}}{\beta(\gamma_n)} \left(
        \frac{\gamma_{n+1} - y \beta(\gamma_n)}{\beta(\gamma_n)-y\gamma_{n+1}}
      \right) .
  \end{alignat}
\end{subequations}
And further define
\begin{align}
  \mathcal{B}_n(a) 
  &= {\X}_n(a) + {\Y}_n(a){\Z}_n(a),
  &\mathcal{C}_n(a) 
  &= {\Y}_n(a) \mathcal{A} _n(a).
\end{align}
Equation~\Ref{eqn hayarec} now becomes:
\begin{equation}
H(\gamma_n (a),y \gamma_n(a)) 
= \mathcal{B}_n + \mathcal{C}_n H (\gamma_{n+1} (a),  y \gamma_{n+1}(a)).
\end{equation}
We obtain a solution for $H(a,ya)$ by iterating the above equation
\begin{equation}
H (a,ya) = \mathcal{B}_0 + \mathcal{C}_0 \mathcal{B}_1
+ \mathcal{C}_0 \mathcal{C}_1 \mathcal{B}_2 + \cdots
= \sum_{n=0}^\infty \mathcal{B}_n(a) \prod_{m=0}^{n-1} \mathcal{C}_m(a).
\label{eqn Haya}
\end{equation}
As was the case for the symmetric wedge, we are able to simplify the
above expression by rewriting $\gamma_n(a)$ in terms of the original
kernel roots and thereby rewrite the expressions for $\mathcal{B}_n$
and $\mathcal{C}_n$.

In order to find $H(a,b)$ from equation~\Ref{eqn Hkern}, we need both
$H(a,ya)$ \emph{and}$H(yb,b)$. Equation~\Ref{eqn halpb} gives $H(b,yb)$ in terms of $H(\alpha_1(b),y \alpha_1(b))$:
\begin{equation}
  H(yb,b) = -\frac{X(\alpha_1(b),b)}{Z(\alpha_1(b),b)} - \frac{Y(\alpha_1(b),b)}{Z(\alpha_1(b),b)} H(\alpha_1(b),y \alpha_1(b))
\end{equation}
So using the above expressions for $H(a,ya)$ and $H(yb,b)$ we have, at
least in principle, a solution for $H(a,b)$:
\begin{multline}
  \label{hab formal}
  H (a,b)  = \frac{X(a,b)}{K(a,b)}   
  -\frac{Z(a,b)  X(\alpha_1(b),b)}{Z(\alpha_1(b),b) K(a,b)} + \frac{Y(a,b)}{K(a,b)} H(a,ya) \\
  - \frac{Z(a,b) Y(\alpha_1(b),b)}{Z(\alpha_1(b),b) K(a,b)} H(\alpha_1(b),y \alpha_1(b))
\end{multline}
Of course, we would like to be able to simplify the above
expression. In particular we would like to rewrite $\gamma_n(a)$ and
$\gamma_n(\alpha_1(b))$ in terms of simpler functions, as we did for
$\beta_n(a)$ in Section~\ref{ssec explicit sym}.

Interestingly enough, it is possible to determine $H(a,ya)$ in
equation~\Ref{eqn Haya} by inspection of the terms in this
expression.  Putting $x=y=t$, one obtains 
\begin{eqnarray}
\hspace{-1.5cm}
H (a,ta) 
& = & - \sum_{n=0}^\infty \left[ \frac{t^{2(n+1)^2-3}}{a} \right]
\left[ \frac{a - \beta t - a \beta t^2 + a\beta t^{2n+2}
}{a(1+\beta)t^{2n}-\beta(a+t^{2n-1})}
\right] \nonumber \\
& & \times \left[
\frac{a-\beta t - a\beta t^2 - \beta t^{4n+1} +a(1+\beta)t^{4n+2}
}{a+\left[\frac{1-t^{2n}}{1-t^2}\right]
(a(1-\beta)t^2+a(1+\beta)t^{2n+2})
-\left[\frac{1-t^{4n}}{1-t^2}\right]\beta t}
\right] \nonumber \\
& & \times
\prod_{m=0}^n
\left[
\frac{a(1+\beta)t^{2m}-\beta(a+t^{2m-1})
}{a-\beta t - a\beta t^2 + a\beta t^{2m+2}}
\right] .
\label{eqn 87} 
\end{eqnarray}
From this expression one may determine $H(bt,t)$ from
equation~\Ref{eqn halpb}, and thus an expression for $H(a,b)$.
While the resulting expression gives a series expansion
for the numbers of paths, it is not very useful because it
is so complex.  In the next section we proceed by 
simplifying expressions for the compositions of
the $\alpha$'s and $\beta$'s above; ultimately
this will lead to a simpler expression for $H(a,b)$.

\subsection{Simplifying things}
In much the same way as for the symmetric case, we can find simple
expressions for the $1/\gamma_n(a)$ in terms of the original kernel
roots. Consideration of the kernel and its roots gives:
\begin{subequations}
  \begin{align}
    \frac{1}{\alpha_{-1}(b)} + \frac{1}{\alpha_{+1}(b)} & = 
    \frac{1+y^2}{yb}, \\
    \frac{1}{\beta_{-1} (a)} + \frac{1}{\beta_{+1} (a)} & = 
    \frac{1+y^2}{ya} - \frac{x(1-y^2)}{y}.
  \end{align}
\end{subequations}
Since certain compositions of $\alpha$ and $\beta$ give the identity (see
equations~\Ref{eqn alpbet comp}), we have the additional relations:
\begin{subequations}
  \begin{align}
    \frac{1}{\alpha_1 ( \beta_1 (a) )} = \frac{1}{\gamma_1(a)} 
    & =  \frac{1+y^2}{y\beta_1 (a)} - \frac{1}{a} , \label{eqn 57a}  \\
    \frac{1}{\beta_1( \alpha_1 (b))} 
    & =  \frac{1+y^2}{y\alpha_1 (b)} - \frac{1}{b} - \frac{x(1-y^2)}{y}.
    \label{eqn 57b} 
  \end{align}
\end{subequations}
We note that the last term in equation~\Ref{eqn 57b} means that the
resulting expressions for $\gamma_n(a)$ are more complicated than
those for $\beta_n(a)$ for the symmetric case (see equation~\Ref{eqn
  betan}); this in turn leads to a significantly more complicated
solution.

Setting $a = \gamma_{n-1}(a)$ and $b = \beta_1 ( \gamma_{n-1}(a))$ in
the above two equations give:
\begin{subequations}
  \begin{align}
    \frac{1}{\gamma_n (a)} &
    = \frac{1+y^2}{y\beta_1(\gamma_{n-1}(a))} - \frac{1}{\gamma_{n-1}(a)},\\
    \frac{1}{\beta_1 (\gamma_n (a))} & = 
    \frac{1+y^2}{y\gamma_n (a)} 
    - \frac{1}{\beta_1 (\gamma_{n-1}(a))} - \frac{x(1-y^2)}{y}.    
  \end{align}
\end{subequations}
These equations can be solved:
\begin{subequations}
  \begin{align}
    \frac{1}{\gamma_n (a)} 
    & = \frac{1-y^{4n}}{y^{2n-1}(1-y^2)\beta_1(a)}
    - \frac{1-y^{4n-2}}{y^{2n-2}(1-y^2)a}
    - \frac{x(1-y^{2n})(1-y^{2n-2})}{y^{2n-2}(1-y^2)}, \nn
    & = \frac{1}{1-y^2}\left(x(1+y^2) + Q(a)\, y^{2n} + y \bar{Q}(a)\, y^{-2n} \right),\\
    \frac{1}{\beta_1(\gamma_n(a))} 
     &= \frac{1-y^{4n+2}}{y^{2n}(1-y^2)\beta_1(a)}
     - \frac{1-y^{4n}}{y^{2n-1}(1-y^2)a}
     - \frac{x(1-y^{2n})^2}{y^{2n-1}(1-y^2)},\nn
    & = \frac{1}{1-y^2}\left(2xy + y Q(a)\,y^{2n} + \bar{Q}(a)\,y^{-2n} \right).
  \end{align}
\end{subequations}
where we have used
\begin{equation}
  Q(a) = \frac{1}{a} - \frac{y}{\beta_1(a)} - x
  \qquad 
  \bar{Q}(a)  = \frac{1}{\beta_1(a)} - \frac{y}{a} - xy
\end{equation}
Note that  $\bar{Q}(a) Q(a) = x^2y$. In fact we can reduce the above
expressions for $\gamma_n$ and $\beta(\gamma_n)$ even further using
this fact:
\begin{subequations}
  \begin{align}
    \frac{1}{\gamma_n(a)} 
    & = \frac{(x+y^{2n-2}Q)(x+y^{2n}Q)}{y^{2n-2} (1-y^2) Q}\\
    \frac{1}{\beta(\gamma_n(a))} 
    & = \frac{(x+y^{2n} Q )^2}{y^{2n-1} (1-y^2) Q}
  \end{align}
\end{subequations}

The above then lead to the following expressions that will be
useful in writing down our solution:
\begin{subequations}
  \begin{align}
    \frac{1}{\gamma_n(a)} - \frac{y}{\beta_1(\gamma_n(a))}
    &= (x + y^{2n}Q), \\
    \frac{1}{\beta_1(\gamma_n(a))} - \frac{y}{\gamma_n(a)}
    &= \frac{x}{y^{2n-1}Q} \left(x + y^{2n} Q \right), \\
    \frac{1}{\beta_1(\gamma_n(a))} - \frac{y}{\gamma_{n+1}(a)}
    &= y \left(x+y^{2n}Q \right), \\
    \frac{1}{\gamma_{n+1}(a)} - \frac{y}{\beta_1(\gamma_n(a))} 
    &= \frac{x}{y^{2n} Q} \left(x + y^{2n} Q \right)
  \end{align}
\end{subequations}
This in turn lets us write
\begin{subequations}
  \begin{align}
    \frac{\beta(\gamma_n) - y \gamma_n}{\gamma_n - y \beta(\gamma_n)} 
    & = \frac{y^{2n-1}}{x} Q \\
    \frac{\gamma_{n+1} - y \beta(\gamma_n)}{\beta(\gamma_n) - y \gamma_{n+1}}
    & = \frac{y^{2n+1}}{x}Q\\
  \end{align}
\end{subequations}
where we have made use of the fact that $\bar{Q} = x^2y / Q$.
Hence $\mathcal{C}_n = \mathcal{Y}_n \mathcal{A}_n$ can now be written as
\begin{equation}
  \mathcal{C}_n = \frac{\gamma_{n+1}}{\gamma_n}
  \left(\frac{\beta(\gamma_n) - y \gamma_n}{\gamma_n - y \beta(\gamma_n)} \right)
  \left(\frac{\gamma_{n+1} - y \beta(\gamma_n)}{\beta(\gamma_n) - y \gamma_{n+1}} \right)
  = \frac{\gamma_{n+1}} {\gamma_n} \cdot \frac{ y^{4n}}{x^2} Q^2
\end{equation}
In a similar way we find that
\begin{subequations}
  \begin{align}
     \frac{\beta(\gamma_n) - y \gamma_n}{\gamma_n^2}
     & = y (x + y^{2n-2} Q)  \\
     \frac{\gamma_{n+1} - y \beta(\gamma_n)}{\gamma_n \gamma_{n+1}} 
     & = y^2(x + y^{2n-2} Q )   
  \end{align}
\end{subequations}
which allows us to also simplify the expression for $\mathcal{B}_n =
\mathcal{X}_n + \mathcal{Y}_n \mathcal{Z}_n$:
\begin{subequations}
  \begin{align}
    \mathcal{X}_n 
    & = \frac{1}{xy} \left(\frac{\beta(\gamma_n) - y \gamma_n}{\gamma_n^2}\right)
    &&= \frac{x + y^{2n-2}Q}{x}   \\
    \mathcal{Y}_n\mathcal{Z}_n 
    & = - \frac{1}{xy} 
    \left( \frac{\gamma_{n+1} - y\beta(\gamma_n) }{\gamma_n
        \gamma_{n+1}} \right)
    \left( \frac{\beta(\gamma_n) - y \gamma_n}{\gamma_n - y
        \beta(\gamma_n)} \right)
    && = -\frac{y^{2n}Q}{x^2} \left( x + y^{2n-2} Q \right)\\
    \mathcal{B}_n 
    & = \mathcal{X}_n + \mathcal{Y}_n \mathcal{Z}_n
    & & = \frac{1}{x^2} \left(x + y^{2n-2} Q \right)\left(x - y^{2n} Q \right)
  \end{align}
\end{subequations}
Substituting the above into the expression for $H(a,ya)$ in
equation~\Ref{eqn Haya} gives:
\begin{align}
  H(a,ya) 
  & = \sum_{n=0}^\infty \mathcal{B}_n(a) \prod_{m=0}^{n-1} \mathcal{C}_m(a).\nn
  & = \sum_{n=0}^\infty 
  \frac{1}{x^2}
  \left(x + y^{2n-2} Q \right)\left(x - y^{2n} Q \right)
  \frac{\gamma_n(a)}{a} \left(\frac{Q}{x}\right)^{2n} y^{2n(n-1)}\nn
  & = \frac{(1-y^2) Q }{a x^2 y^2} \sum_{n=0}^\infty 
  \frac{(x-y^{2n}Q) }{(x+y^{2n}Q)}
  \left(\frac{Q}{x} \right)^{2n}
  y^{2n^2}
\end{align}
which is a significant simplification of equation~\Ref{eqn 87}.
We can now substitute this into equation~\Ref{hab formal} to obtain
$H(a,b)$. This requires us to compute 
$H(\alpha_1(b), y \alpha_1(b) )$ from the above expression. 
Let
\begin{equation}
  P(b) = Q(\alpha_1(b)) = 
\frac{y}{2b}\left( 1-2xyb-y^2+\sqrt{(1-y^2)(1-4xyb-y^2)}\right)
\end{equation}
 then we have
\begin{equation}
  H(\alpha, y \alpha) = \frac{(1-y^2) P }{\alpha_1(b) x^2 y^2} \sum_{n=0}^\infty 
  \frac{(x-y^{2n}P) }{(x+y^{2n}P)}
  \left(\frac{P}{x} \right)^{2n}
  y^{2n^2}
\end{equation}

Below we give the length generating function (when
$x=t,y=t,a=1,b=1$).
\begin{prop}
  The generating function of partially directed walks ending in a
  horizontal step in the wedge $\mathcal{W}_1$ is
  \begin{align}
    \label{eqn h11 full}
    h_1 (1,1) 
    &= \frac{(1-t)^2 - \sqrt{(1-t^2)(1-5t^2)}}{2(1-2t-t^2)} \nn
    & - Q \frac{(1-t^2) }{t^2(1-2t-t^2)} 
    \sum_{n=0}^\infty 
    \frac{(1-t^{2n-1}Q) }{(1+t^{2n-1}Q)}
    \left(\frac{Q}{t} \right)^{2n}
    t^{2n^2}\nn
    & + \frac{(1-t^2)}{1-2t-t^2}
    \sum_{n=0}^\infty 
    \frac{(1-t^{2n-1}P) }{(1+t^{2n-1}P)}
    \left(\frac{P}{t} \right)^{2n}
    t^{2n^2}
  \end{align}
  where $t$ counts the number of edges and 
  \begin{subequations}
    \begin{align}
      Q(1;t,t) & =(1-t-t^2-t^3-\sqrt{(1-t^4)(1-2t-t^2)})/2 \\
      P(1;t,t) &= (1-3t^2-\sqrt{(1-t^2)(1-5t^2)})/2t
    \end{align}
  \end{subequations}
  The generating function of all paths in $\mathcal{W}_1$ is then
  $k_1(1,1;t,t) = (h_1(1,1;t,t)-1)/t$.
\end{prop}

As was the case for the symmetric wedge, the functions $P$ and 
$Q$ that make up our expression for $H$ have combinatorial 
interpretations in terms of partially directed walks bounded 
by a single line.  Let $B_-(x,y)$ be the generating function
of walks that end with a horizontal step, start and end on 
the line $Y=0$ and stay on or above that same line. Then
\begin{equation}
  Q(a;x,y) = xy^2 (B_-(ax,y) - 1) .
\end{equation}
Similarly let $B_/(x,y)$ be the generating function of walks 
that end with a horizontal step, start and end on
the line $Y=X$ and stay on or above that same line. Then
\begin{equation}
  P(b;x,y) = xy^2 (B_/(bx,y)  - 1) .
\end{equation}
Again we would like to find a more direct combinatorial 
derivation of the generating functions $H(1,1)$ and $H(1,t)$. 
We have been unable to do so.

\subsection{Asymptotics for \texorpdfstring{$p=1$}{p=1}} 

Before we study the asymptotics of walks in the wedge $\mathcal{W}_1$,
let us compute the number of walks lying on or above the line $y=0$.
\begin{lemma}
  \label{lem upper bound}
  The generating function of partially directed walks lying on or
  above the line $y=0$ is
  \begin{equation}
    \frac{-1+z+3z^2+z^3-\sqrt{(1-z^4)(1-2z-z^2)}}{2z^2(z^2-2z-1)}
  \end{equation}
  The number of these walks is asymptotic to
  \begin{equation}
    \sqrt{ \frac{7+5\sqrt{2}}{2\pi}} \frac{ (1+\sqrt{2})^n }{\sqrt{n}}
    \left( 1 + O(1/n) \right)
  \end{equation}
  Hence the number of walks in the wedge $\mathcal{W}_1$ must also be
  $O( (\sqrt{2}+1)^n / \sqrt{n} )$.
\end{lemma}
\proof One can derive a functional equation for the generating
function of such walks which can be solved using the kernel method.
The asymptotics can then be computed by analysing the dominant
singularity at $t=\sqrt{2}-1$. The last result follows since the
number of walks in the wedge $\mathcal{W}_1$ cannot exceed the number
of walks lying above $y=0$. \qed

As was the case for walks in the symmetric wedge, we analyse the
asymptotics of partially directed walks in the asymmetric wedge
$\mathcal{W}_1$ by singularity analysis. Let us split the expression
given in equation~\Ref{eqn h11 full} into 3 pieces and study their
dominant singularities:
\begin{subequations}
  \begin{align}
    p_1 &= \frac{(1-t)^2 - \sqrt{(1-t^2)(1-5t^2)}}{2(1-2t-t^2)} \\
    p_2 &= - Q \frac{(1-t^2) }{t^2(1-2t-t^2)} \sum_{n=0}^\infty
    \frac{(1-t^{2n-1}Q) }{(1+t^{2n-1}Q)} \left(\frac{Q}{t}
    \right)^{2n} t^{2n^2} \\
    p_3 &= \frac{(1-t^2)}{1-2t-t^2} \sum_{n=0}^\infty \frac{(1-t^{2n-1}P)
    }{(1+t^{2n-1}P)} \left(\frac{P}{t} \right)^{2n} t^{2n^2}
  \end{align}
\end{subequations}

Let us treat the asymptotics of each of these functions separately.
\begin{lemma}
  The coefficients of $p_1(t)$ are asymptotic to
  \begin{equation}
    [t^n] p_1 = -\sqrt{\frac{5}{8 \pi}} \cdot 
    \left( (2+\sqrt{5}) - (-1)^n(\sqrt{5}-2) \right)
    \cdot
    \frac{ \left(\sqrt{5} \right)^n}{\sqrt{n^3}} 
    \cdot
    \left(1 +  O(n^{-1}) \right)
  \end{equation}
\end{lemma}
\proof The generating function $p_1$ appears to have 6 singularities:
2 simple poles from the zeros of the denominators and four square-root
singularities at $t=\pm 1, \pm 1/\sqrt{5}$. Closer analysis shows that
there are no singularities at the zeros of the denominator and that
generating function is dominated by the singularities at $t = \pm
1/\sqrt{5}$. Analysis (by the techniques in \cite{Flajolet1990} ) of
these singularities leads to the above expression. \qed

Before we can study the asymptotics of $p_2(t)$ and $p_3(t)$ we need
the following lemma about the location of the zeros of $1+Qt^k$ and $1+Pt^k$.
\begin{lemma}
  For $k=-1,0,1,2,\ldots$, the functions $1+Q(t)t^k$ and $1+P(t)t^k$ are
  not zero within the disk $|t|<1/2$.
\end{lemma}
\proof The function $Q(t)$ satisfies $ Q^2 - (1-t-t^2-t^3) Q +t^4 = 0
$, and so $h = Q t^k$  satisfies:
\begin{equation}
  h^2 - (1-t-t^2-t^3) ht^k + t^{2k+4}.
\end{equation}
Now if $1+Qt^k=0$ we have $h=-1$ and so
\begin{equation}
  \label{eqn Q zero}
  1 + (1-t-t^2-t^3) t^k + t^{2k+4} = 0.
\end{equation}
For $|t| \leq 1/2$ we can bound $|1-t-t^2-t^3| \leq (1+\frac{1}{2} +
\frac{1}{4} + \frac{1}{8}) < 2$, and therefore
\begin{equation}
  | (1-t-t^2-t^3) t^k + t^{2k+4} | \leq |2 t^k + t^{2k+4} | \leq 3
  |t|^k < 3 \cdot 2^{-k}
\end{equation}
Hence for $k \geq 2$, the above quantity is less than $1$ and so
equation~\Ref{eqn Q zero} cannot be satisfied. It remains to check the
cases $k=-1,0,1$. In these cases we can solve equation~\Ref{eqn Q
  zero} directly and verify that $t$ lies outside $|t|\leq 1/2$.

The argument for $P(t)$ follows a similar line. The function $h = P t^k$ satisfies:
\begin{equation}
  h^2 - t^{k-1}(1-3t^2) h + t^{2k+2} = 0.
\end{equation}
Hence if $h=-1$ we have
\begin{equation}
  \label{eqn P zero}
  1 + (1-3t^2) t^{k-1} + t^{2k+2} = 0
\end{equation}
For $|t|<1/2$ we can bound $|1-3t^2| < 2$ and so
\begin{equation}
  | (1-3t^2) t^{k-1} + t^{2k+2} | \leq  3 |t|^{k-1}
\end{equation}
Hence for $k\geq3$, equation~\Ref{eqn P zero} cannot be satisfied for
$|t| \leq 1/2$. For $k=-1,0,1,2$, we can check
equation~\Ref{eqn P zero} directly and verify that the zeros do not
lie inside $|t|<1/2$. 

Note that when $k=0$, equation~\Ref{eqn P zero} has a solution
$t=\sqrt{2}-1$, however this point corresponds to the other branch of
$P$ being $+1$.

 \qed

We can now move onto the asymptotics of $p_2(t)$ and $p_3(t)$.
\begin{lemma}
  \label{lem p3 asympt}
  The coefficients of $p_3(t)$ are asymptotic to
  \begin{align}
    [t^n] p_3(t) & = \frac{(\sqrt{2}+1)^n}{\sqrt{2}}
    \sum_{k=0}^\infty
    \frac{1-(\sqrt{2}-1)^{2k+1}}{1+(\sqrt{2}-1)^{2k+1}}
    (\sqrt{2}-1)^{2k^2+2k} + O\left( \left(\sqrt{5} \right)^n \right) \nn
    & = (0.31096381899209832\ldots )(\sqrt{2}+1)^n + O\left( \left(\sqrt{5}
    \right)^n \right)
  \end{align}
\end{lemma}
\proof The function $p_3$ has a simple pole from the zero of the
denominator of prefactor. There are also square-root singularities in
$P(t)$ at $t = \pm 1, \pm 1/\sqrt{5}$, simple poles when $1+P
t^{2n-1}=0$ and a natural boundary at $|t|=1$ coming from the
$\theta$-function like structure of the sum. Of these singularities,
the simple pole dominates, followed by the singularities at $\pm
1/\sqrt{5}$. The simple pole and its residue give the dominant
asymptotics and the square-root singularities give the $O(5^{n/2})$
corrections. \qed

\begin{lemma}
  \label{lem p2 asympt}
  The coefficients of $p_2(t)$ are asymptotic to
  \begin{align}
    [t^n] p_2(t) & = -\frac{(\sqrt{2}+1)^n}{\sqrt{2}}
    \sum_{k=0}^\infty
    \frac{1-(\sqrt{2}-1)^{2k+1}}{1+(\sqrt{2}-1)^{2k+1}}
    (\sqrt{2}-1)^{2k^2+2k} \left(1 + o(1) \right) \nn
    & = -(0.31096381899209832\ldots )(\sqrt{2}+1)^n \left(1 + o(1) \right)
  \end{align}
  The dominant term is equal in magnitude, but opposite in sign, to
  the dominant term in the asymptotics of $p_3(t)$.
\end{lemma}
\proof The singularities of the function $p_2(t)$ arise from
singularities of the prefactor, the singularities of $Q(t)$, the zeros
of $1+Q t^{2n-1}$ and the natural boundary at $|t|=1$ from the
$\theta$-function structure of the sum. Of these, the simple pole of
the prefactor at $t=\sqrt{2}-1$ and the singularities arising from
$Q(t)$ at the same point, dominate the asymptotics.

The contribution from the simple pole may be computed by finding its
residue. We note that $P(\sqrt{2}-1) = Q(\sqrt{2}-1) = 3-2\sqrt{2}$,
and this means that the residue is in fact equal in magnitude, but
opposite in sign, to that computed for $p_3(t)$.

\qed

We note that if the dominant asymptotics of $p_2(t)$ and $p_3(t)$ must
cancel each other. Otherwise, the number of walks in this wedge is
$\sim (\sqrt{2}-1)^n$ which would contradict Lemma~\ref{lem upper bound}.

\begin{lemma}
  The $k^{\mathrm{th}}$ summand of $p_2(t)$ is
  \begin{equation}
    p_{2,k}(t) = \left( - Q \frac{(1-t^2) }{t^2(1-2t-t^2)} 
       \frac{(1-t^{2k-1}Q) }{(1+t^{2k-1}Q)} \left(\frac{Q}{t}
       \right)^{2k} t^{2k^2} \right) 
  \end{equation}
  so that $p_2(t) = \sum p_{2,k}(t)$. The coefficient of $t^n$ in
  $p_{2,k}(t)$ is asymptotic to
  \begin{multline}
    [t^n] p_{2,k}(t) = 
    - (1+\sqrt{2})^{n} \cdot \frac{1}{\sqrt{2}}
    \left( \frac{1-(\sqrt{2}-1)^{2k+1}}{1+(\sqrt{2}-1)^{2k+1}} \right)
    ( 1+\sqrt{2} )^{-2k^2-2k} \\
    +
    (1+\sqrt{2})^n \cdot 
    \sqrt{\frac{2}{\pi n}} 
    \left[
    \frac{(2k+1)(1-(\sqrt{2}-1)^{4k+2})-(\sqrt{2}-1)^{2k+1}}{
        (1+(\sqrt{2}-1)^{2k+1})^2} \right.\\
    \left.
      + O\left(\frac{1}{\sqrt{n}}\right)
    \right] (1+\sqrt{2})^{-2k^2-2k-5/2}
  \end{multline}
\end{lemma}
\proof The result follows from standard singularity analysis
(\cite{Flajolet1990}) of $p_{2,k}(t)$. \qed

\begin{rem}
  It is unfortunately the case that we have been unable to proceed
  completely rigorously from this point. In particular, we have been
  unable to obtain uniform bounds on the error terms in the above
  asymptotic expressions. On the basis of numerical testing, we think
  that the that the error term is $O(k^3/\sqrt{n})$. If this is the
  case, then one can sum the contributions of the individual
  $p_{2,k}(t)$ to obtain the asymptotics of coefficients of $p_2(t)$.
 
  We believe that the expressions that follow are indeed \emph{exact},
  if not completely rigorous.
\end{rem}

Assuming that the asymptotic expression in the previous
lemma has a uniform error bound, so that we may sum the contribution
to the individual $p_{2_k}(t)$, we find that
\begin{multline} 
  [t^n] p_2(t) =  (1 + \sqrt{2})^n \Big(-0.31096381899209832\ldots \\
  + \frac{0.090584741026764287\ldots}{\sqrt{n}}  + O\big(1/\sqrt{n^3}\big) \Big)
\end{multline}
where the constant $0.31\dots$ is the constant that appears in both
Lemmas~\ref{lem p3 asympt} and~\ref{lem p2 asympt}. So adding together
the contributions from the $p_i$ we obtain
\begin{equation}
  [t^n] h_1(1,1) = (1 + \sqrt{2})^n \left( \frac{0.090584741026764287\ldots}{\sqrt{n}}  + O\big(1/\sqrt{n^3}\big)  \right)
\end{equation}
The generating function $h_1(1,1)$ enumerates walks that end in a
horizontal step, and so we obtain the number of walks in
$\mathcal{W}_1$ ending with any step by multiplying this expression by
a factor of $(1+\sqrt(2))$ (since the generating functions differ by a
factor of $t$):
\begin{equation}
  [t^n] k_1(1,1) = (1 + \sqrt{2})^n \left( \frac{0.218693916694303177\ldots}{\sqrt{n}}  + O\big(1/\sqrt{n^3}\big)  \right)
\end{equation}
We have confirmed this numerically using the first 1000 terms in the
series expansion of $k_1(1,1)$. Additionally we have check the above
expression using Bruno Salvy's \emph{gdev} package for Maple
\cite{Salvy91b}.

\begin{rem}
  We note that if the above result is indeed made rigorous then we
  have a result analogous to Corollary~\ref{cor symm bound}.  The
  number of walks in any wedge $\mathcal{W}_p$ for $1 \leq p < \infty$
  is bounded between the number of walks in $\mathcal{W}_1$ and
  partially directed walks inside the first quadrant (see
  Lemma~\ref{lem upper bound}. Hence the number of partially directed
  walks in the wedge $\mathcal{W}_p$, $c_n^{(p)}$ obeys the following
  inequality
  \begin{equation}
    0.21869\ldots \leq \lim_{n \to \infty} \frac{c_n^{(p)} n^{1/2} }{(1+\sqrt{2})^n} \leq
     \sqrt{\frac{7+5\sqrt{2}}{2\pi}} = 1.496489\ldots
  \end{equation}
  for any $1 \leq p < \infty$.
\end{rem}

\section{Conclusions}
In this paper we have proved that partially directed paths in the
wedges $\mathcal{V}_p$ and $\mathcal{W}_p$ all grow with the same
exponential growth rate $1 + \sqrt{2}$ independent of $p$.
Additionally we have found generating functions for partially directed
paths in the symmetric wedge $\mathcal{V}_1$ and the asymmetric wedge
$\mathcal{W}_1$, using a variation of the kernel method. From these
generating functions we have computed the asymptotics of the number of
paths in both of these wedges.

Curiously the number of paths in the symmetric wedge, $\mathcal{V}_1$,
has the same leading asymptotic behaviour as partially directed paths
with no bounding wedge. Similarly the number of paths in the
asymmetric wedge, $\mathcal{W}_1$, has the same leading asymptotic
behaviour as partially directed paths above the line $Y=0$. Because of
this, we are able to determine the leading asymptotic behaviour of
paths in the wedges $\mathcal{V}_p$ and $\mathcal{W}_p$ for all $p
\geq 1$.


\noindent\section*{Acknowledgements}
EJJvR is supported by a grant from National Sciences and Engineering
Research Council of Canada. AR is partially supported by funding from
the Australian Reserach Council and the Centre of Excellence for
Mathematics and Statistics of Complex Systems.

\bibliographystyle{plain}
\bibliography{wedge_ref2.bib}

\begin{thebibliography}{10}

\bibitem{Andre1887}
D.~{Andr\'e}.
\newblock Solution directe du {probl\`eme r\'esolu} par {M. Bertrand}.
\newblock {\em Comptes Rendus Acad. Sci. Paris}, 105:436--437, 1887.

\bibitem{B02}
C.~Banderier, M.~Bousquet-M\'elou, A.~Denise, P.~Flajolet, D.~Gardy, and
  D.~Gouyou-Beauchamps.
\newblock Generating functions for generating trees.
\newblock {\em Discrete Math.}, 246:29--55, 2002.

\bibitem{BM96}
M.~Bousquet-M\'elou.
\newblock A method for the enumeration of various classes of column-convex
  polygons.
\newblock {\em Discrete Math.}, 154:1--25, 1996.

\bibitem{BM2001}
M.~Bousquet-M\'elou.
\newblock Walks on the slit plane: other approaches.
\newblock {\em Adv. in Appl. Math.}, 27:243--288, 2001.

\bibitem{BM02}
M.~Bousquet-M\'elou.
\newblock Counting walks in the quarter plane.
\newblock In {\em Mathematics and Computer Science: Algorithms, trees,
  combinatorics and probabilities}, Trends in Mathematics, pages 49--67.
  Birkhauser, 2002.

\bibitem{BM03}
M.~Bousquet-M\'elou.
\newblock Four classes of pattern-avoiding permutations under one roof:
  generating trees with two labels.
\newblock {\em Electronic J. Combinatorics}, 9(2):R19, 2003.

\bibitem{BM2002a}
M.~Bousquet-M\'elou and G.~Schaeffer.
\newblock Walks on the slit plane.
\newblock {\em Probab. Theory Relat. Fields}, 124:305--344, 2002.

\bibitem{Cardy1984}
J.L Cardy.
\newblock Conformal invariance and surface critical behavior.
\newblock {\em Nuc. Phys. B.}, 240:514--532, 1984.

\bibitem{DiMarzio1971}
E.A. DiMarzio and R.J. Rubin.
\newblock Adsorption of a chain polymer between two plates.
\newblock {\em J. Chem. Phys.}, 55:4318--4336, 1971.

\bibitem{D00}
P.~Duchon.
\newblock On the enumeration of generation of generalized dyck words.
\newblock {\em Discrete Math.}, 225:121--135, 2000.

\bibitem{Duplantier1986}
B.~Duplantier and H.~Saleur.
\newblock Exact surface and wedge exponents for polymers in two dimensions.
\newblock {\em Phys. Rev. Lett.}, 57:3179--3182, 1986.

\bibitem{Barcucci2001}
R.~Pinzani E.~Barcucci, E.~Pergola and S.~Rinaldi.
\newblock A bijection for some paths on the slit plane.
\newblock {\em Adv. in Appl. Math.}, 26:89--96, 2001.

\bibitem{Fayolle1999}
G.~Fayolle, R.~Iasnogorodski, and V.~Malyshev.
\newblock {\em Random walks in the quarter plane: algebraic methods, boundary
  value problems and applications}.
\newblock Number~40 in Applications of Mathematics. Springer Verlag, Berlin,
  1999.

\bibitem{Flajolet1990}
P.~Flajolet and A.~M. Odlyzko.
\newblock Singularity analysis of generating functions.
\newblock {\em SIAM J. Discrete Math.}, 3:216--240, 1990.

\bibitem{Gessel1992}
I.~M. Gessel and D.~Zeilberger.
\newblock Random walk in a weyl chamber.
\newblock {\em Proc. Amer. Math. Soc.}, 115:27--31, 1992.

\bibitem{HW85}
J.~M. Hammersley, , and S.~G. Whittington.
\newblock Self-avoiding walks in wedges.
\newblock {\em J. Phys. A: Math. Gen.}, 18:101--111, 1985.

\bibitem{K63}
H.~Kesten.
\newblock On the number of self-avoiding walks.
\newblock {\em J. Math. Phys.}, 4:960--959, 1963.

\bibitem{K64}
H.~Kesten.
\newblock On the number of self-avoiding walks ii.
\newblock {\em J. Math. Phys.}, 5:1128--1137, 1964.

\bibitem{MS93}
N.~Madras and G.~Slade.
\newblock {\em The Self-Avoiding Walk}.
\newblock Birkhauser, London, 1993.

\bibitem{MW40}
W.~K. McCrae and F.~J.~W. Whipple.
\newblock Random paths in two and three dimensions.
\newblock {\em Proc. Roy. Soc. Edinburgh}, 60:281--298, 1940.

\bibitem{N83}
D.~H. Napper.
\newblock {\em Polymeric Stabilization of Colloidal Dispersions}.
\newblock Academic Press, London, 1983.

\bibitem{Niederhausen1998}
H.~Niederhausen.
\newblock Lattice paths between diagonal boundaries.
\newblock {\em Electronic J. Combinatorics}, 5:R30, 1998.

\bibitem{Brak2005}
A.~Rechnitzer R.~Brak, A.L.~Owczarek and S.G. Whittington.
\newblock A directed walk model of a long chain polymer in a slit with
  attractive walls.
\newblock {\em J. Phys. A: Math. Gen.}, 38:4309--4325, 2005.

\bibitem{Salvy91b}
B.~Salvy.
\newblock Examples of automatic asymptotic expansions.
\newblock {\em SIGSAM Bulletin}, 25(2):4--17, April 1991.

\bibitem{JvR05c}
E.~J.~Janse van Rensburg.
\newblock Forces in motzkin paths in a wedge.
\newblock To appear in J. Phys. A: Math. Gen.

\bibitem{JvR05b}
E.~J.~Janse van Rensburg.
\newblock Adsorbing bargraph paths in a $q$-wedge.
\newblock {\em J. Phys. A: Math. Gen.}, 38:8505--8525, 2005.

\bibitem{JvR05a}
E.~J.~Janse van Rensburg.
\newblock Square lattice directed paths adsorbing on the line $y=qx$.
\newblock {\em J. Stat. Mech.: Theo. and Exp.}, page P09010, 2005.

\bibitem{JvRL05}
E.~J.~Janse van Rensburg and Y.~Le.
\newblock Forces in square lattice directed paths in a wedge.
\newblock {\em J. Phys. A: Math. Gen.}, 38:8493--8503, 2005.

\bibitem{V98}
C.~Vanderzande.
\newblock {\em Lattice Models of Polymers}, volume~11 of {\em Cambride Lecture
  Notes in Physics}.
\newblock Cambridge University Press, 1998.

\end{thebibliography}

\end{document}